\documentclass[11pt]{article}
\usepackage{amsmath,amsfonts,amsthm,amssymb,ifthen}
\usepackage{fancybox,color,epsfig}

\newboolean{@details}
\setboolean{@details}{true}

\def\detail#1{\ifthenelse{\boolean{@details}}{\marginpar{begin\\
details}#1\marginpar{end\\
details}}{}}

\newboolean{@todoon}
\setboolean{@todoon}{true}

\def\todo#1{\ifthenelse{\boolean{@todoon}}{\textbf{(\color{red} todo\footnote{\color{red} #1})}}{}}

\newtheorem{theorem}{Theorem}[section]

\newtheorem{lemma}[theorem]{Lemma}

\newtheorem{proposition}[theorem]{Proposition}
\newtheorem{definition}[theorem]{Definition}

\newcommand{\sdd}{$\mathrm{SDD}_{0}$}
\newcommand{\sddmz}{$\mathrm{SDDM}_{0}$}
\newcommand{\sddm}{$\mathrm{SDDM}$}

\newenvironment{tightlist}%
{\begin{list}{}{\usecounter{bean}%
  \setlength{\topsep}{0em}
  \setlength{\parsep}{0em}
  \setlength{\partopsep}{0em}
  \setlength{\itemsep}{0em}
}}%
{\end{list}}

\newenvironment{dasenumerate}%
{
\begin{enumerate}}%
{\end{enumerate}}

\def\pleq{\preccurlyeq}

\newcommand\uu{\boldsymbol{\mathit{u}}}
\newcommand\vv{\boldsymbol{\mathit{v}}}

\newcommand\xx{\boldsymbol{\mathit{x}}}

\def\xxtil{\boldsymbol{\mathit{\tilde{x}}}}

\renewcommand\ss{\boldsymbol{\mathit{s}}}
\newcommand\sss{\boldsymbol{\mathit{s}}}
\newcommand\yy{\boldsymbol{\mathit{y}}}
\newcommand\zz{\boldsymbol{\mathit{z}}}
\newcommand\rr{\boldsymbol{\mathit{r}}}
\newcommand\bb{\boldsymbol{\mathit{b}}}

\def\Ahat{\widehat{A}}

\def\bbhat{\hat{\bb}}
\def\sshat{\hat{\ss}}
\def\xxhat{\hat{\xx}}
\def\yyhat{\hat{\yy}}
\def\zzhat{\hat{\zz}}

\def\Bhat{\widehat{B}}
\def\Lhat{\widehat{L}}

\def\prob#1#2{\Pr_{#1}\left[ #2 \right]}

\def\norm#1{\left\| #1 \right\|}

\def\infnorm#1{\left\| #1 \right\|_{\infty }}

\def\setof#1{\left\{#1  \right\}}
\def\sizeof#1{\left|#1  \right|}

\def\bvec#1{{\mbox{\boldmath $#1$}}}

\def\pleq{\preccurlyeq}

\def\abs#1{\left|#1  \right|}
\def\intersect{\cap}

\def\solve#1#2{\mathtt{solve}_{#1} (#2)}
\def\Solve#1{\mathtt{solve}_{#1}}

\newcommand{\ceiling}[1]{\left\lceil#1\right\rceil}

\def\bigO#1{O\left(#1  \right)}

\def\setof#1{\left\{#1  \right\}}
\def\abs#1{\left|#1  \right|}

\newcommand\xxt{\boldsymbol{\mathit{\tilde{x}}}}

\newdimen\pIR
\newcommand\StevesR{{\rm I\kern\pIR R}}
\def\Reals#1{\StevesR^{#1}}

\def\pinv#1{{#1}^{\dagger}}
\def\inv#1{{#1}^{-1}}

\newenvironment{fminipage}%
  {\begin{Sbox}\begin{minipage}}%
  {\end{minipage}\end{Sbox}\fbox{\TheSbox}}

\newenvironment{algbox}[0]{\vskip 0.2in
\noindent 
\begin{fminipage}{6.3in}
}{
\end{fminipage}
\vskip 0.2in
}

\def\edg#1{\pmb{\boldsymbol{(}} #1 \pmb{\boldsymbol{)}}}

\def\noff#1{\mathrm{noff}\left(#1 \right)}
\def\dim#1{\mathrm{dim}\left(#1 \right)}

\def\stretch#1#2{\textrm{st}_{#1} (#2)}
\def\weight#1{\textrm{weight} (#1)}
\def\res#1{\textrm{resistance} (#1)}

\def\union{\cup}
\def\intersect{\cap}

\begin{document}

\title{
Nearly-Linear Time Algorithms for Preconditioning and Solving
Symmetric, Diagonally Dominant Linear Systems%
\thanks{
This paper is the last in a sequence of three papers expanding
  on material that appeared first under the title
  ``Nearly-linear time algorithms for graph partitioning, 
    graph sparsification, and solving linear systems''~\cite{SpielmanTengPrecon}.
The second paper, ``Spectral Sparsification of Graphs''~\cite{SpielmanTengSparsifier}
  contains algorithms for constructing sparsifiers of graphs, which we
  use in this paper to build preconditioners.
The first paper,
``A Local Clustering Algorithm for Massive Graphs and its Application to Nearly-Linear Time Graph Partitioning''~\cite{SpielmanTengCuts}
 contains graph partitioning algorithms that are used to construct sparsifiers in the
  second paper.
\vskip 0.01in
This material is based upon work supported by the National Science Foundation 
  under Grant Nos. 0325630, 0324914, 0634957, 0635102, 0707522, 0964481, 1111257 and 1111270.
Any opinions, findings, and conclusions or recommendations expressed in this material are those of the authors and do not necessarily reflect the views of the National Science Foundation.
}
}

\author{
Daniel A. Spielman \\ 
Department of Computer Science\\
Program in Applied Mathematics\\
Yale University
\and 
Shang-Hua Teng \\
Department of Computer Science\\
Viterbi School of Engineering\\
University of Southern California}

\maketitle

\begin{abstract}
We present a randomized algorithm that, 
  on input a symmetric, weakly diagonally dominant 
  $n$-by-$n$ matrix $A$ with $m$ nonzero entries
  and an $n$-vector $\bb$, produces an $\xxtil $ such that
  $\norm{\xxtil - \pinv{A} \bb }_{A} \leq \epsilon \norm{\pinv{A} \bb }_{A}$
  in expected time
  $O ( m \log^{c}n \log (1/\epsilon)),$
  for some constant $c$.
By applying this algorithm inside the inverse power method, we compute
  approximate Fiedler vectors in a similar amount of time.
The algorithm
  applies subgraph preconditioners in a recursive fashion.
These preconditioners improve upon the subgraph preconditioners
  first introduced by Vaidya (1990).

For any symmetric, weakly diagonally-dominant matrix $A$ 
  with non-positive off-diagonal entries and $k \geq 1$, we construct
  in time $O (m \log^{c} n)$
 a preconditioner $B$ of $A$
  with at most 
$2 (n - 1) + O ((m/k) \log^{39} n)$
nonzero off-diagonal entries such that
  the finite generalized condition number $\kappa_{f} (A,B)$
  is at most $k$,
for some other constant $c$.

In the special case when the nonzero structure of the matrix is planar
 the corresponding linear system solver runs in expected time
$  O (n \log^{2} n + n \log n \  \log \log n \ \log (1/\epsilon ))$.

We hope that our introduction of algorithms of low asymptotic complexity
  will lead to the development of algorithms that are also fast 
  in practice.
\end{abstract}

\section{Introduction}\label{sec:intro}
We
  design an algorithm with nearly optimal asymptotic complexity 
  for solving linear systems in
  symmetric, weakly diagonally dominant (\sdd) matrices.
The algorithm applies a classical iterative
  solver, such as the Preconditioned Conjugate Gradient or the
  Preconditioned Chebyshev Method,
  with a novel preconditioner that we construct
  and analyze using techniques from graph theory.
Linear systems in these preconditioners may be reduced to systems
  of smaller size in linear time by use of a direct method.
The smaller linear systems are solved recursively.
The resulting algorithm solves linear systems in \sdd \
  matrices in time that is asymptotically almost linear in their number of
  nonzero entries.
Our analysis does not make any assumptions about the nonzero structure
  of the matrix, and thus may be applied to the solution of the systems
  in  \sdd \ matrices that arise in any application,
  such as the solution of elliptic partial differential 
  equations by the finite element method~\cite{StrangITAM,BomanHendricksonVavasis},
  the solution of maximum flow problems by interior point 
  algorithms~\cite{FrangioniGentile,DaitchSpielman}, or the solution of 
  learning problems on graphs~\cite{BelkinRegression,ZhouBLWS03,ZhuGL}.

\begin{figure}[h]
\begin{center}
\hspace{0.5in}
\begin{minipage}{.3\linewidth}
$\begin{bmatrix}
1.5 & - 1.5 & 0 & 0\\
-1.5 & 4 & -2 & -0.5\\
0 & -2 & 3 & -1\\
0 & -0.5 & -1 & 1.5
\end{bmatrix}
$
\end{minipage}\hspace{1in}
\begin{minipage}{.3\linewidth}
\epsfig{figure=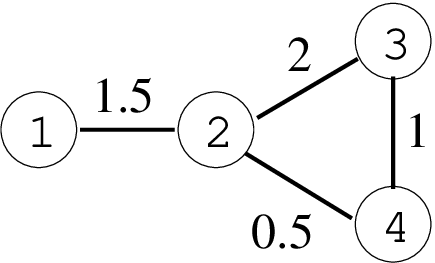,height=1in}
\end{minipage}
\end{center}
\caption{A Laplacian matrix and its corresponding weighted graph.}
\label{fig:aGraph}
\end{figure}

Graph theory drives the construction of our preconditioners.
Our algorithm is best understood by first examining its behavior
  on Laplacian matrices---symmetric matrices
  with non-positive off-diagonals and zero row sums.
Each $n$-by-$n$ Laplacian matrix $A$ may be associated with a weighted graph,
  in which the weight of the edge between distinct
  vertices $i$ and $j$ is $-A_{i,j}$ (see Figure~\ref{fig:aGraph}).
We precondition the Laplacian matrix $A$ of a graph $G$ by the Laplacian matrix $B$
  of a subgraph $H$ of $G$ that resembles a spanning tree of $G$ plus a few edges.
The subgraph $H$ is called an \textit{ultra-sparsifier} of $G$, and its corresponding
  Laplacian matrix is a very good preconditioner for $A$:
The finite generalized condition number $\kappa_{f} (A,B)$ is $\log^{O (1)} n$.
Moreover, it is easy to solve linear equations in $B$.
As the graph $H$ resembles a tree plus a few
  edges, we may use partial Cholesky factorization to eliminate most of the
  rows and columns of $B$ while incurring only a linear amount of fill.
We then solve the reduced system recursively.

This paper contains two principal contributions.
The first, which appears in Sections~\ref{sec:solvers} through~\ref{sec:stability},
  is the analysis of a multilevel algorithm that uses ultra-sparsifiers to solve
  systems of equations in \sdd \ matrices.
The second, which appears in Sections~\ref{sec:graphs}
  through~\ref{sec:ultra}, is the construction of ultra-sparsifiers
  for Laplacian matrices.
In the remainder of the introduction we formally define ultra-sparsifiers
  and the sparsifiers from which they are built.
In Section~\ref{sec:prior}, we survey the contributions upon which we
  build, the improvements upon our work, and other related work.
We devote Section~\ref{sec:background} to 
  recalling the basics of support theory,
  explain
  how the problem of solving linear equations in \sdd \ matrices may be reduced
  to that of solving equations in positive definite \sdd \ matrices with nonnegative
  off-diagonal entries, and explain how the problem of preconditioning such matrices
  can be reduced to that of preconditioning
  Laplacian matrices.

In Section~\ref{sec:solvers}, we state the properties we require of
  partial Cholesky factorizations, and we present our first algorithms for
  solving equations in \sdd-matrices.
These algorithms directly solve equations in the preconditioners, rather
  than using a recursive approach, and take time
  roughly $O (m^{5/4} \log^{c} n)$, for some constant $c$, for general 
  \sdd-matrices and time $O (n^{9/8} \log^{1/2}  n)$ for 
  \sddm-matrices with planar nonzero structure.
To accelerate these algorithms, we apply our preconditioners in a recursive fashion.
We analyze the complexity of these recursive algorithms in Section~\ref{sec:alg}, 
  obtaining our main algorithmic results.
In Section~\ref{sec:stability}, we show that our algorithm provides
  accurate answers even with computations in limited precision.
In Section~\ref{sec:eigs}, we observe that these linear system solvers yield
  efficient algorithms for computing
  approximate Fiedler vectors, when applied inside the inverse power method.

\textit{We do not attempt to optimize the exponent of $\log n$ in the complexity
  of our algorithm.}
Rather, we present the simplest analysis we can find of an algorithm of complexity
  $$O ( m \log^{c}n \log (1/\epsilon))$$ for some constant $c$.
Recently, Koutis, Miller and Peng~\cite{KMP1,KMP2} have discovered much simpler
  constructions of ultra-sparsifiers which lead to algorithms achieving
  such a running time for every $c > 1$.
Similarly, we do not prove tight bounds on the precision required for our algorithms to work.
We merely prove that $O (\log \kappa (A) \log^{c} n \log \epsilon^{-1})$ bits of precision suffice.
We hope that the efficient preconditioners of Koutis, Miller and Peng~\cite{KMP1,KMP2}
  will motivate a tighter stability analysis.


\subsection{Definitions and Notation}
We recall that a matrix $A$ is \textit{weakly diagonally dominant} if
  $A (i,i) \geq \sum_{j \not = i} \abs{A (i,j)}$ for all $i$.
We define \sdd \ to be the class of symmetric, weakly diagonally dominant
  matrices, and \sddmz \  to be the class of \sdd-matrices with non-positive off-diagonal
  entries.
We let \sddm \ be the class of positive-definite \sddmz \ matrices.
The \sddm-matrices are M-matrices and in particular are
  Stieltjes matrices.
A Laplacian matrix is a \sddmz-matrix with zero row-sums.

Throughout this paper, we define the $A$-norm by
\[
  \norm{\xx}_{A} = \sqrt{\xx^{T} A \xx}.
\]
For symmetric matrices $A$ and $B$, we write
\[
A \pleq B
\]
if $B - A$ is positive semidefinite.
We recall that if $A$ is positive semidefinite and $B$ is symmetric, then
  all eigenvalues of $A B$ are real.
For a matrix $B$, we let $\pinv{B}$ denote the Moore-Penrose
  pseudo-inverse of $B$---that is the matrix with the same nullspace
  as $B$ that acts as the inverse of $B$ on its image.
We let $\kappa (A)$ denote the ratio of the largest to the smallest singular
  value of $A$.
These are the largest and smallest eigenvalues when $A$ is symmetric
  and positive definite.

We denote the logarithm base $2$ of $x$ by $\log x$ and the natural
  logarithm of $x$ by $\ln  x$.

\subsection{Ultra-sparsifiers}

\begin{definition}[Ultra-Sparsifiers]\label{def:ultraSparsifier}
A $(k,h)$-\textit{ultra-sparsifier} of an 
  $n$-by-$n$ 
  \sddm-matrix $A$ with 
  $2m$ nonzero off-diagonal entries  is 
  a  \sddm-matrix $A_{s}$
  such that
\begin{dasenumerate}
\item $  A_{s} \pleq A  \pleq  k \cdot A_{s}$.
\item $A_{s}$ has at most $2 (n-1) + 2 h m/k$ nonzero off-diagonal entries.
\item The set of  nonzero entries of
  $A_{s}$ is a subset of the set of nonzero
  entries of $A$.
\end{dasenumerate}
\end{definition}

In Section~\ref{sec:ultra}, we present a randomized algorithm that runs in
   expected time $O (m \log^{c} n)$ 
  that takes as input a Laplacian matrix
  $A$ and a $k \geq 1$ and produces a $(k,h)$-ultra-sparsifier
  of $A$ with probability at least $1 - 1/2n$, for
\begin{equation}\label{eqn:h}
 h = c_{3} \log_{2}^{c_{4}} n,
\end{equation}
where $c$, $c_{3}$ and $c_{4}$ are some absolute constants.
As we will use these ultra-sparsifiers throughout the paper, we will define
  a $k$-\textit{ultra-sparsifier} to be a $(k,h)$-ultra-sparsifier
  where $h$ satisfies \eqref{eqn:h}.

For matrices whose graphs are planar, we 
  present a simpler construction of $(k,h)$-ultra-sparsifiers,
  with $h = \bigO{\log n \ \log \log n}$.
This simple constructions exploits low-stretch spanning
trees~\cite{AKPW,EEST,AbrahamBartalNeiman,AbrahamNeiman}, and is presented in
  Section~\ref{sec:lowStretch}.
Our construction of ultra-sparsifiers in Section~\ref{sec:ultra} builds
  upon the simpler construction, but requires the use of
  \textit{spectral sparsifiers} \cite{SpielmanTengSparsifier}.
The following definition of sparsifiers will suffice for the purposes
  of this paper.

\begin{definition}[Spectral Sparsifiers]\label{def:sparsifier}
A $d$-\textit{sparsifier} of  an $n$-by-$n$  \sddm-matrix $A$ is 
  a  \sddm-matrix $A_{s}$
  such that
\begin{dasenumerate}
\item  $  A_{s} \pleq A \pleq (5/4) A_{s}$.
\item  $A_{s}$ has at most $d n$ nonzero
  off-diagonal entries.
\item  The set of  nonzero entries of
  $A_{s}$ is a subset of the set of nonzero
  entries of $A$.
\item  For all $i$,
\[
  \sum_{j \not = i} 
\frac{A_{s} (i,j)}{A (i,j)}
\leq 
2 \sizeof{\setof{j : A (i,j) \not = 0}}.
\]
\end{dasenumerate}
\end{definition}

In a companion paper~\cite{SpielmanTengSparsifier},
  we present a randomized algorithm \texttt{Sparsify2}
  that produces sparsifiers of Laplacian matrices in expected 
  nearly-linear time.
As explained in Section~\ref{sec:background},
  this construction can trivially be extended to all \sddm-matrices.
\begin{theorem}[Spectral Sparsification]\label{thm:sparsifier}
On input an $n\times n$ Laplacian matrix $A$ with $2m$ nonzero
   off-diagonal entries and a $p > 0$,
  \texttt{Sparsify2} runs 
  in expected time $O (m \log (1/p) \log^{17} n)$
  and
  with probability at least $1-p$
  produces a $c_{1} \log^{c_{2}} (n/p)$-sparsifier of $A$,
for $c_{2} = 34$ and some
  absolute constant $c_{1} > 1$.
\end{theorem}


Spielman and Srivastava~\cite{SpielmanSrivastava}
  construct sparsifiers with $c_{2} = 1$, but their construction requires
  the solution of linear equations in Laplacian matrices and
  so can not be used to help speed up the algorithms in this paper.
Their algorithm can be made reasonably fast by using the linear systems
  solvers of Koutis, Miller and Peng~\cite{KMP2}.
Batson, Spielman and Srivastava~\cite{BatsonSpielmanSrivastava} have
  proved that there exist sparsifiers that satisfy conditions $(a)$ through $(c)$ of
  Definition~\ref{def:sparsifier} with $c_{2} = 0$.  

\section{Related Work}\label{sec:prior}

In this section, we explain how our results relate to other
  rigorous asymptotic analyses of algorithms for solving systems
  of linear equations.
For the most part, we restrict our attention to algorithms that
  make structural assumptions about their input matrices, rather
  than assumptions about the origins of those matrices.

Throughout our discussion, we consider an $n$-by-$n$ matrix with
  $m$ nonzero entries.
When $m$ is large relative to $n$ and the matrix is arbitrary,
  the fastest algorithms for solving linear equations are those based on
  fast matrix multiplication~\cite{CoppersmithWinograd,Stothers,Vassilevska}, which take time
 approximately $O (n^{2.37})$.
The fastest algorithm for solving general sparse
  positive semidefinite linear systems
  is the Conjugate Gradient.
Assuming computation with infinite precision,
  one can show that it obtains the correct answer after
  $O (m n)$ 
  operations (see ~\cite[Theorem~38.3]{TrefethenBau}).
To the best of our knowledge, every faster algorithm requires
  additional properties of the input matrix.

\subsection{Special nonzero structure}
In the design and analysis of direct solvers, it is standard
  to represent the nonzero structure of a matrix $A$
  by an unweighted graph $G_{A}$
  that has an edge between vertices $i \neq j$ if and only if
  $A_{i,j}$ is nonzero (see~\cite{DuffErismanReid}).
If this graph has special structure,
  there may be elimination orderings that accelerate direct solvers.
If $A$ is tri-diagonal, in which case $G_{A}$ is a path,
  then a linear system in $A$ can be solved in time $O (n)$.
Similarly, when $G_{A}$ is a tree a linear system in $A$ can be solved in time $O (n)$
  (see \cite{DuffErismanReid}).
 
If the graph of nonzero entries $G_{A}$ is planar,
  one can use Generalized Nested Dissection~\cite{George,LiptonRoseTarjan,GilbertTarjan}
  to find an elimination ordering under which
  Cholesky factorization can be performed in time
  $O (n^{1.5})$ and produces factors with at most $O (n \log n)$ nonzero entries.
We will exploit these results in our algorithms for solving planar
  linear systems in Section~\ref{sec:solvers}.
We recall that a planar graph on $n$ vertices has at most $3n - 6$ 
 edges (see \cite[Corollary 11.1 (c)]{Harary}), so $m \leq 6n$.

\subsection{Subgraph Preconditioners}
Our work builds on a remarkable approach to solving linear systems in
  Laplacian matrices introduced by Vaidya~\cite{Vaidya}.
Vaidya
  demonstrated that a good
  preconditioner for a Laplacian matrix
  $A$ 
  can be found in the Laplacian matrix $B$ of a subgraph of the graph corresponding
  to $A$.
He then showed that one could bound the condition number
  of the preconditioned system by bounding the dilation and congestion
  of an embedding of the graph of $A$ into the
  graph of $B$. 
By using preconditioners 
  obtained by adding edges to maximum spanning trees,
  Vaidya developed
  an algorithm that finds $\epsilon$-approximate solutions to 
  linear systems in \sddmz-matrices
  with at most $d$ nonzero entries per row in 
  time $O ((d n)^{1.75} \log (1 / \epsilon ))$.
For matrices whose corresponding graphs have special
  structure,
  such as having a
  bounded genus or avoiding certain minors,
  he obtained even faster algorithms.
For example, his algorithm for solving planar systems
  runs in time $O ((d n)^{1.2} \log (1 / \epsilon ))$.

As Vaidya's paper was never published
  and his manuscript lacked many
  proofs, the task of formally working out his results fell to others.
Much of its
  content appears in the thesis of his student, Anil Joshi~\cite{Joshi}.
Chen and Toledo~\cite{ChenToledo} present an experimental study
  of  Vaidya's preconditioners,
  and a complete exposition of Vaidya's work along with many extensions was presented
  by Bern \textit{et. al.}~\cite{SupportGraph}.
Gremban, Miller  and Zagha~\cite{Gremban,GrembanMillerZagha}
  explain parts of Vaidya's paper as well
  as extend Vaidya's techniques.
Among other results, they find ways of constructing preconditioners by
  \textit{adding} vertices to the graphs.
Maggs \textit{et. al.}~\cite{MaggsEtAl} prove that
  this technique may be used to construct
  excellent preconditioners,
  but it is still not clear if they can be constructed efficiently.

Gremban~\cite[Lemma 7.3]{Gremban} (see also Appendix~\ref{sec:gremban}) 
  presents a reduction from the problem of solving linear systems in
  \sdd  \ matrices to that of solving linear systems in \sddmz \ matrices
  that are twice as large.
The machinery needed to apply Vaidya's techniques directly
  to matrices with positive off-diagonal elements is developed
  in~\cite{MWB}.
An algebraic extension of Vaidya's techniques for bounding
  the condition number was presented 
  by Boman and Hendrickson~\cite{SupportTheory},
  and later used by them~\cite{BomanHendricksonAKPW}
  to prove that the low-stretch spanning trees constructed by 
  Alon, Karp, Peleg, and West \cite{AKPW},
  yield preconditioners 
  for which the preconditioned system has condition number at most
  $m 2^{\bigO{\sqrt{\log n\log\log n}}}$.
They thereby obtained a solver for
  \sddmz \ linear
   systems that produces $\epsilon $-approximate solutions in
  time $m^{1.5 + o (1)} \log (1 / \epsilon )$.
Through improvements in the construction of low-stretch spanning 
  trees~\cite{EEST,AbrahamBartalNeiman,AbrahamNeiman} and a careful analysis
  of the eigenvalue distribution of the preconditioned system,
  Spielman and Woo~\cite{SpielmanWoo} show that when the Preconditioned Conjugate
  Gradient is applied with the best low-stretch spanning tree preconditioners,
  the resulting linear system solver takes time at most 
  $O(m n^{1/3} \log^{1/2} n \log (1/\epsilon))$.
The preconditioners in the present paper are formed by adding edges
  to these low-stretch spanning trees.

The recursive application of subgraph preconditioners
  was pioneered in the work of Joshi~\cite{Joshi} and Reif~\cite{Reif}.
Reif~\cite{Reif} showed how to recursively apply Vaidya's preconditioners
  to solve linear systems in
   \sddmz-matrices with 
  planar nonzero structure and
  at most a constant number of nonzeros
  per row in time
  $O (n^{1 + \beta } \log^{c} (\kappa (A) / \epsilon ))$,
  for some constant $c$,
  for every $\beta > 0$.
While Joshi's analysis is numerically much cleaner, he only analyzes preconditioners
  for simple model problems.
Our recursive scheme uses ideas from both these works, with some simplification.
Koutis and Miller~\cite{KoutisMiller}
  have developed recursive algorithms that solve
  linear systems in  \sddmz-matrices with planar nonzero structure in
  time $O (n \log (1/\epsilon))$.

Koutis, Miller and Peng~\cite{KMP1,KMP2} have recently made substantial
 improvements in the construction of ultra-sparsifiers that result
  in algorithms for solving linear equations in \sdd \ matrices
  in time $O (m \log n \log^{2} \log n \log (1/\epsilon))$.
Their construction has the added advantage of being much simpler than ours.
Slightly better constructions of ultra-sparsifiers have been shown to exist
  by Kolla, Makarychev, Saberi, and Teng~\cite{KMST}, although their construction
  takes longer than nearly-linear time.

\subsection{Other families of matrices}

Subgraph preconditioners have been used to solve systems of linear equations
  from a few other families.

Daitch and Spielman~\cite{DaitchSpielman} have shown how to reduce
  the problem of solving linear equations in symmetric $M_{0}$-matrices
  to the problem of solving linear equations in \sddmz-matrices,
  given a factorization of the $M_{0}$-matrix of width 2~\cite{FactorWidth}.
These matrices, with the required factorizations, arise in the solution of
  the generalized maximum flow problem by interior point algorithms.

Shklarski and Toledo~\cite{ShklarskiToledo} introduce an extension of
  support graph preconditioners, called \textit{fretsaw preconditioners},
  which are well suited to preconditioning finite element matrices.
Daitch and Spielman~\cite{DaitchSpielmanTruss} use these preconditioners
  to solve linear equations in the stiffness matrices of two-dimensional truss
  structures in time $O (n^{5/4} \log n \log (1/\epsilon))$.

For linear equations that arise when solving
  elliptic partial differential equations, other techniques
  supply fast algorithms.
For example, Multigrid methods 
  provably run in nearly-linear time
 when applied to the solution of some of these
  linear systems \cite{MultigridBook},
  and algorithms based on $\mathcal{H}$-matrices run in nearly-linear time when
  given a sufficiently nice discretization~\cite{Hierarchical}.
Boman, Hendrickson, and Vavasis 
  \cite{BomanHendricksonVavasis}
  have shown that the problem of solving 
  a large class of these linear systems 
  may be reduced to that of solving
  diagonally-dominant systems.
Thus, our algorithms may be applied to the solution of these systems.

\section{Background}\label{sec:background}

We will use the following propositions, whose proofs are elementary.

\begin{proposition}\label{pro:pleqNull}
If $A$ and $B$ are positive semidefinite matrices such that for some
  $\alpha, \beta > 0$,
\[
  \alpha  A \pleq B \pleq \beta A
\]
then $A$ and $B$ have the same nullspace.
\end{proposition}

\begin{proposition}\label{pro:pleqPinv}
If $A$ and $B$ are positive semidefinite matrices 
  having the same nullspace and $\alpha  > 0$, then
\[
    \alpha A \pleq B
\]
if and only if
\[
  \alpha \pinv{B} \pleq \pinv{A}.
\]
\end{proposition}

The following proposition establishes the equivalence of two notions
  of preconditioning.
This proposition is called the ``Support Lemma'' in~\cite{SupportGraph}
  and~\cite{Gremban}, and is implied by Theorem~10.1 of~\cite{Axelsson}.

\begin{proposition}\label{pro:pleq}
If $A$ and $B$ are symmetric positive semidefinite matrices
  with the same nullspace,
  then
 all eigenvalues of $A \pinv{B}$
  lie between $\lambda _{min}$ and $\lambda _{max}$
 if and only if
\[
  \lambda _{min} B \pleq A \pleq \lambda _{max} B.
\]
\end{proposition}

Following Bern \textit{et. al.}~\cite{SupportGraph}, we define the
  \textit{finite generalized condition number} $\kappa_{f} (A,B)$
  of matrices $A$ and $B$ having the same nullspace
  to be the ratio of the largest to smallest nonzero eigenvalues
  $A \pinv{B}$.
Proposition~\ref{pro:pleq} tells us that
  $\lambda_{min} B \pleq A \pleq \lambda_{max} B$ implies
  $\kappa_{f} (A,B) \leq \lambda_{max} / \lambda_{min}$.
One can use $\kappa_{f} (A,B)$ to bound the number of iterations
  taken by the Preconditioned Conjugate Gradient algorithm
  to solve linear systems in $A$ when using $B$ as a preconditioner.
Given bounds on $\lambda_{max}$ and $\lambda_{min}$, one can similarly
  bound the complexity of the Preconditioned Chebyshev method.

\subsection{Preconditioning}\label{ssec:laplacians}

When constructing preconditioners, we will focus our attention
  on the problem of preconditioning Laplacian matrices.
Bern \textit{et. al.}~\cite{SupportGraph},
  observe that the problem of preconditioning  \sddmz-matrices
  is easily reduced to that of preconditioning Laplacian matrices.
We recall the reduction as we will make use of it later.

Any \sddmz -matrix $A$ can be decomposed as $A = A_{L} + A_{D}$ where $A_{L}$ is a Laplacian
  matrix and $A_{D}$ is a diagonal matrix with non-negative entries.
Given a Laplacian matrix $B_{L}$ that preconditions $A_{L}$,
  we use $B = B_{L} + A_{D}$ as a preconditioner for $A$.

\begin{proposition}[\cite{SupportGraph}, Lemma 2.5]\label{pro:reduceToLaplacians}
Let $A$ be a  \sddmz-matrix
  and let $A = A_{L} + A_{D}$ where $A_{L}$ is a Laplacian
  matrix and $A_{D}$ is a diagonal matrix with non-negative entries.
If $B_{L}$ is another Laplacian matrix on the same vertex set, then
\[
  \kappa_{f} (A , B_{L} + A_{D}) \leq \kappa_{f} (A_{L}, B_{L}).
\]
In particular,
  if $A_{L} \pleq B_{L}$ then
  $A \pleq B_{L} + A_{D}$.
Similarly, 
  if
  $B_{L} \pleq A_{L}$, then $B_{L} + A_{D} \pleq A$.
\end{proposition}

So, any algorithm for constructing sparsifiers or ultra-sparsifiers
  for Laplacian matrices can immediately be converted into an algorithm
  for constructing sparsifiers or ultra-sparsifiers for \sddmz-matrices.
This is why we 
  restrict our attention to the problem of preconditioning Laplacian matrices
  in Sections~\ref{sec:lowStretch} and~\ref{sec:ultra}.

\subsection{Solving equations}\label{ssec:solving}
In this section, we describe how one can quickly transform the general
  problem of solving a system of equations in a \sdd -matrix to the problem
  of solving a system of equations in an irreducible \sddm -matrix.

Recall that a symmetric matrix $A$ is \textit{reducible} if there is a permutation
  matrix $P$ for which $P^{T} A P$ is a block-diagonal matrix with at least two
  blocks.
If such a permutation exists, one can find it in linear time.
A matrix that is not reducible is said to be \textit{irreducible}.
The problem of solving a linear system in a reducible matrix can be
  reduced to the problems of solving linear systems in each of the blocks.
Accordingly, we will only consider the problem of solving
  linear systems in irreducible matrices.
It is well-known that a symmetric matrix is irreducible if and only if
  its corresponding graph of nonzero entries is connected.
We use this fact in the special case of Laplacian matrices,
  observing that the weighted graph associated with a Laplacian matrix $A$
  has the same set of edges as $G_{A}$.

\begin{proposition}\label{pro:connected}
A Laplacian matrix is irreducible if and only if its
  corresponding weighted graph is connected.
\end{proposition}

It is also well-known that the null-space of the Laplacian
  matrix of a connected graph is the span of the all-1's vector.
Combining this fact with Proposition~\ref{pro:reduceToLaplacians},
  one can show that the only singular irreducible
  \sddm-matrices are the Laplacian matrices.

\begin{proposition}\label{pro:singularM0}
A singular irreducible \sddm-matrix is a Laplacian matrix,
 and its nullspace is spanned by the all-1's vector.
\end{proposition}

We now note that by Gremban's reduction, the problem of solving an
  equation of the form
  $A \xx = \bb $ for a \sdd-matrix $A$ can be reduced
  to the problem of solving a system that is twice as large in a
   \sddm-matrix (see Appendix~\ref{sec:gremban}), without
  any loss of approximation quality.
So, for the purposes of asymptotic complexity
  we need only consider the problem of solving systems in
   \sddm-matrices.

While the algorithms we develop may be naturally applied to the solution of
  equations in both positive definite and singular \sddmz \ matrices,
  it is simpler to analyze the algorithms by considering just one of
  these cases.
We find it simpler to reduce the singular, Laplacian, case to the
  positive definite case, and then to analyze our solvers for
  positive definite matrices.
Let $A$ be an irreducible Laplacian matrix.
As the nullspace of $A$ is spanned by the constant vectors, the
  equation $A \xx = \bb$ will only have a solution if the sum
  of the entries of $\bb$ is zero.
In this case, the system is under-determined and for every solution $\xx$
  the vector $\xx - \bvec{1} \xx (1)$ is a solution as well.
Thus, we may assume that $\xx (1) = 0$ and seek a solution in
  the remaining variables.
If we let $A_{2}$ be the submatrix of $A$ containing all but
  its first row and column and we let $\xx_{2}$ and $\bb_{2}$ denote the 2nd
  through last entry of $\xx$ and $\bb$,
  then we find 
\[
  A_{2} \xx_{2} = \bb_{2}.
\]
It is easy to see that this system is positive-definite:
  $A_{2}$ is a diagonally-dominant \sddm-matrix, and 
 the rows corresponding to vertices that are neighbors of vertex $1$
  are strictly diagonally dominant.
If we obtain an approximate solution to this system
  $\xxt_{2}$ and set $\xxt$ to be zero in its first coordinate
  and $\xxt_{2}$ in the rest, then
\[
  \norm{\xxt - \xx}_{A}
= 
  \norm{\xxt_{2} - \xx_{2}}_{A_{2}}.
\]
As $\norm{\xx}_{A} = \norm{\xx_{2}}_{A_{2}}$, the guarantee that our solver returns
  an $\xxt_{2}$ satisfying
  $  \norm{\xxt_{2} - \xx_{2}}_{A_{2}} \leq \epsilon \norm{\xx_{2}}_{A_{2}}$
  implies that
  $  \norm{\xxt - \xx}_{A} \leq \epsilon \norm{\xx}_{A}$.
As we have assumed that $A$ is irreducible, its nullspace is just the span of the  
  constant vector.
So, we can bring $\xxt$ close to $\pinv{A} \bb$ by subtracting the average entry of $\xxt$
  from each of its entries.
However, this is not strictly necessary as we have
\[
  \norm{\xxt - \pinv{A} \bb}_{A} =  \norm{\xxt - \xx}_{A} \leq \epsilon \norm{\xx}_{A} 
  = \epsilon \norm{\pinv{A} \bb}_{A}.
\]

So, our bound on the quality of $\xxt$ as a solution to the singular system
  is the same as our bound on the quality of $\xxt_{2}$ as a solution to the positive-definite
  system.

\section{Solvers and One-Level Algorithms}\label{sec:solvers}

To solve a system in an irreducible  \sddm-matrix 
  $A$, we will compute an ultra-sparsifier  $B$ of $A$,
  and then solve the system in $A$ using a preconditioned iterative method.
At each iteration of this method, we will need to solve a system in $B$.
We will solve a system in $B$ by a two-step algorithm.
We will first apply Cholesky factorization repeatedly to eliminate all 
  rows and columns with at most one or two nonzero off-diagonal entries.
As we stop the Cholesky factorization before it has factored
  the entire matrix, we call
 this process a \textit{partial Cholesky factorization}.
We then apply another solver on the remaining system.
In this section, we analyze the use of a direct solver.
In Section~\ref{sec:alg}, we obtain our fastest algorithms by solving
  the remaining system recursively.

\subsection{Partial Cholesky Factorization}
The application of partial Cholesky factorization to eliminate rows
  and columns with at most 2 nonzero off-diagonal entries
  results in a factorization of $B$ of the form
\[
  B = P L C L^{T} P^{T},
\]
where $C$ has the form
\[
  C = \left(
  \begin{array}{ll}
 I_{n-n_{1}} & 0 \\
  0 & A_{1},
\end{array}
 \right),
\]
$P$ is a permutation matrix, $L$ is nonsingular and lower triangular of the form
\[
  L = \left(
  \begin{array}{ll}
 L_{1,1} & 0 \\
 L_{2,1} & I_{n_{1}},
\end{array}
 \right),
\]
and every row and column of $A_{1}$ has
  at least $3$ nonzero off-diagonal entries.

In the following proposition we state properties of this factorization that we will exploit.
Variants of this proposition are implicit in earlier work on subgraph preconditioners~\cite{Vaidya,Joshi,SupportGraph}.

\begin{proposition}[Partial Cholesky Factorization]\label{pro:chol}
If $B$ is an irreducible  \sddm-matrix then,
\begin{dasenumerate}
\item  $A_{1}$ is an irreducible  \sddm-matrix.

\item If the graph of nonzero entries of $B$ is planar,
  then the graph of nonzero entries of $A_{1}$ is as well.

\item $L$ has at most $3n$ nonzero entries.

\item If $B$ has $2 (n - 1 + j)$ nonzero off-diagonal entries, then
  $A_{1}$ has dimension at most $2 j - 2$ and has at most $2 (3 j - 3)$ 
  nonzero off-diagonal entries.
\end{dasenumerate}
\end{proposition}
\begin{proof}
It is routine to verify that $A_{1}$ is diagonally dominant with non-positive
  off-diagonal entries, and that planarity is preserved by elimination of
  rows and columns with $2$ or $3$ nonzero entries, as these correspond
  to vertices of degree $1$ or $2$ in the graph of nonzero entries.
It is similarly routine to observe that these eliminations preserve irreducibility
  and singularity.

To bound the number of entries in $L$, we note that for each row and column with
  $1$ nonzero off-diagonal entry
  that is eliminated, the corresponding column in $L$
  has 2 nonzero entries, 
  and that for each row and column with $2$ nonzero off-diagonal entries
  that is eliminated, the corresponding column in $L$
  has 3 nonzero entries.

To bound $n_{1}$, the dimension of $A_{1}$, first observe that the elimination of a 
  row and column with $1$ or $2$ nonzero off-diagonal entries
  decreases both the dimension by $1$ and the number of nonzero entries
  by $2$.
So, $A_{1}$ will have $2 (n_{1} -1 + j)$ nonzero off-diagonal entries.
As each row in $A_{1}$ has at least $3$ nonzero off-diagonal entries,
  we have
\[
    2 (n_{1} - 1 + j) \geq 3 n_{1},
\]
which implies $n_{1} \leq 2j - 2$.
The bound on the number nonzero off-diagonal entries in $A_{1}$ 
  follows immediately.

\end{proof}

We name the algorithm that performs this factorization
  $\mathtt{PartialChol}$, and invoke it with the syntax
\[
  (P, L, A_{1}) = \mathtt{PartialChol} (B).
\]
We remark that \texttt{PartialChol} can be implemented to run
  in linear time.

\subsection{One-Level Algorithms}\label{sec:onelevel}

Before analyzing the algorithm in which we solve systems in $A_{1}$
  recursively, 
  we pause to examine the complexity of an algorithm that
  applies a direct solver to systems in $A_{1}$.
While the results in this subsection are not necessary for the main
  claims of our paper, we hope they will provide intuition.

If we are willing to ignore numerical issues, we may apply
  the conjugate gradient algorithm to directly solve systems
  in $A_{1}$ in $O (n_{1} m_{1})$ operations~\cite[Theorem~38.3]{TrefethenBau}, 
  where $m_{1}$ is the number
  of nonzero entries in $A_{1}$.
In the following theorem, we examine the performance of the resulting algorithm.

\begin{theorem}[General One-Level Algorithm]\label{thm:oneShot}
Let $A$ be an irreducible $n$-by-$n$  \sddm-matrix with $2 m$
  nonzero off-diagonal entries.
Let $B$ be a $\sqrt{m}$-ultra-sparsifier of $A$.
Let $(P, L, A_{1}) = \mathtt{PartialChol} (B)$.
Consider the algorithm that 
  solves systems in $A$ by applying PCG with $B$
  as a preconditioner, and solves each system in $B$
  by a performing backward substitution on its partial Cholesky factor,
  solving the inner system in $A_{1}$
  by conjugate gradient used as an exact solver,
  and performing forward substitution on its partial Cholesky factor.
Then for every right-hand side $\bb$, after 
\[
O (m^{1/4} \log (1/\epsilon))
\]
  iterations, comprising
\[
O (m^{5/4} \log^{2 c_{4}} n \log (1/\epsilon))
\]
  arithmetic operations,
the algorithm will output 
 an approximate solution $\xxt$ 
  satisfying
\begin{equation}\label{eqn:oneShot}
\norm{\xxt-\inv{A} \bb}_{A} \leq  \epsilon \norm{\inv{A} \bb}_{A}.
\end{equation}
\end{theorem}
\begin{proof}
As $\kappa_{f} (A,B) \leq \sqrt{m}$, we may apply the
  standard analysis of PCG~\cite{Axelsson}, to show that 
  \eqref{eqn:oneShot} will be satisfied after
  $O (m^{1/4} \log (1/\epsilon))$
  iterations.
To bound the number of operations in each iteration,
  note that $B$ has at most $2 (n-1)+O (\sqrt{m} \log^{c_{4}} n)$
  nonzero off-diagonal entries.
So, Proposition~\ref{pro:chol} implies $m_{1}$ and $n_{1}$ are both
  $O (\sqrt{m} \log^{c_{4}} n)$.
Thus, the time required to solve each inner system in $A_{1}$
  by the Conjugate Gradient is at most 
  $O (m_{1} n_{1}) = O (m \log^{2 c_{4}} n)$.
As $A$ is irreducible, $m \geq n -1$, and so
  this upper bounds the number of operations that must be performed in
  each iteration.
\end{proof}

When the graph of nonzero entries of $A$ is planar, we may
  precondition  using 
  the algorithm \texttt{UltraSimple}, 
  presented in Section~\ref{sec:lowStretch},   
  instead of \texttt{UltraSparsify}.
As the matrix $A_{1}$ produced by applying partial Cholesky
  factorization to the output of \texttt{UltraSimple}
  is also planar, we can solve the linear systems in $A_{1}$
  by the
  generalized nested dissection algorithm of Lipton, Rose
  and Tarjan~\cite{LiptonRoseTarjan}.
This algorithm uses graph separators to choose a good order
  for Cholesky factorization.
The Cholesky factorization is then computed in time
  $O (n_{1}^{3/2})$.
The resulting Cholesky factors only have 
  $O (n_{1} \log n_{1})$ nonzero entries, and so each
  linear system in $A_{1}$ may be solved in time
  $O (n_{1} \log n_{1})$, after the Cholesky factors have been
  computed.

\begin{theorem}[Planar One-Level Algorithm]\label{thm:planarOneShot}
Let $A$ be an $n$-by-$n$ planar, irreducible
   \sddm-matrix
  with $2 m$ nonzero off-diagonal entries.
Consider the algorithm that solves linear systems in $A$
  by using PCG with the preconditioner 
\[
B = \mathtt{UltraSimple} (A, n^{3/4} \log^{1/3} n),
\]
solves systems in $B$ by 
  applying {$\mathtt{PartialChol}$} to 
  factor $B$ into $P L [I, 0; 0, A_{1}] L^{T} P^{T}$,
 and
  uses generalized nested dissection to solve systems in $A_{1}$.
For every right-hand side $\bb$, 
  this algorithm computes an $\xxt$ satisfying
\begin{equation}\label{eqn:planarOneShot}
\norm{\xxt-\pinv{A} \bb}_{A} \leq  \epsilon \norm{\pinv{A} \bb}_{A}
\end{equation}
  in time
\[
  O \left(n^{9/8} \log^{1/2} n \log (1/\epsilon) \right).
\]
\end{theorem}
\begin{proof}
First, recall that the planarity of $A$ implies $m \leq 3n$.
Thus, the time taken by \texttt{UltraSimple} is dominated by the time taken
  by \texttt{LowStretch}, which is $O (n \log n \ \log \log n)$ (see Theorem~\ref{thm:lowStretch}).

By Theorem~\ref{thm:lowStretch} and Theorem~\ref{thm:augmentTree}, the matrix
  $B$ has at most $2 (n-1) + 6 n^{3/4} \log^{1/3} n$ nonzero off-diagonal entries and
\[
\kappa_{f} (A,B)
= 
  O \left(n^{1/4} \log^{2/3} n \log \log n\right) 
  \leq 
  O \left(n^{1/4} \log n\right) .
\]
Again, standard analysis of PCG~\cite{Axelsson}
  tells us that
  the algorithm will require at most
\[
  O \left(n^{1/8} \log^{1/2} n \log (1/\epsilon )\right)
\]
iterations to guarantee that \eqref{eqn:planarOneShot} is satisfied.

By Proposition~\ref{pro:chol}, 
  the dimension of $A_{1}$, $n_{1}$, is at most $6 n^{3/4} \log^{1/3} n$.
Before beginning to solve the linear system, the algorithm will
  spend 
\[
 O (n_{1}^{3/2}) = 
 O ((n^{3/4} \log^{1/3} n)^{3/2}) = O (n^{9/8} \log^{1/2} n)
\]
time
  using generalized nested dissection~\cite{LiptonRoseTarjan}
  to permute and Cholesky factor
  the matrix $A_{1}$.
As the factors obtained will have at most
  $O (n_{1} \log n_{1})  \leq O (n)$ nonzeros,
  each iteration of the PCG will require at most
  $O (n)$ steps.
So, the total complexity of the application of the PCG will be
\[
O \left(n \cdot \left(n^{1/8} \log^{1/2} n \log (1/\epsilon )\right)   \right)
=
  O \left(n^{9/8} \log^{1/2} n  \log (1/\epsilon) \right),
\]
which dominates the time required to compute the Cholesky factors and
  the time of the call to \texttt{UltraSimple}.
\end{proof}
 
We remark that the algorithm of Lipton, Rose and Tarjan~\cite{LiptonRoseTarjan}
  can be accelerated by the use of algorithms for fast matrix inversion~\cite{CoppersmithWinograd,Stothers,Vassilevska}.
One can similarly accelerate our planar one-level algorithm.

\section{The Recursive Solver}\label{sec:alg}

In our recursive algorithm for solving linear equations,
  we solve linear equations in a matrix $A$ by computing an ultra-sparsifier $B$,
  using partial Cholesky factorization to reduce it to a matrix $A_{1}$,
  and then solving the system in $A_{1}$ recursively.
Of course, we compute all of the necessary ultra-sparsifiers and Cholesky factorizations
  just once at the beginning of the algorithm.

In this section we assume infinite precision arithmetic.
We defer an analysis of the impact of limited precision to the next section.

To specify the recursive algorithm for an $n$-by-$n$ matrix,
  we first set the parameter
\begin{equation}\label{eqn:k}
k = (14 h + 1)^{2},
\end{equation}
where we recall that the parameter $h$
  is determined by the quality of the ultra-sparsifiers we can compute
  (see equation \eqref{eqn:h}),

We use the algorithm \texttt{BuildPreconditioners}
  to build the sequence of preconditioners
  and Cholesky factors.
In Section~\ref{sec:ultra}, we define the routine \texttt{UltraSparsify}
  for weighted graphs, and thus implicitly for Laplacian matrices.
To define \texttt{UltraSparsify} 
  for general irreducible \sddm-matrices $A$,
  we express $A$ as a sum of matrices $A_{L}$ and $A_{D}$
  as explained in Proposition~\ref{pro:reduceToLaplacians},
  and set
\[
\mathtt{UltraSparsify} (A, k) =   A_{D} + \mathtt{UltraSparsify} (A_{L}, k).
\]

\begin{algbox}
\noindent \texttt{BuildPreconditioners}$(A_{0})$,
\begin{enumerate}
\item  Set $i = 0$, $h = c_{3} \log_{2}^{c_{4}}\dim{A_{0}}$ (as in \eqref{eqn:h}) and $k = (14 h+1)^{2}$ (as in \eqref{eqn:k}).

\item Repeat
\label{alg:build:until}

\begin{enumerate}
\item $i = i + 1$.
\item  $B_{i} = \mathtt{UltraSparsify} (A_{i-1}, k)$.
\item  $(P_{i}, L_{i}, A_{i}) = \mathtt{partialChol} (B_{i})$.

\end{enumerate}

\item []Until $A_{i}$ has dimension less than $66 h  + 6$.

\item Set $\ell = i$.
\item Compute $Z_{\ell} = \inv{A_{\ell }}$.
\end{enumerate}
\end{algbox}

We now make a few observations about the sequence of matrices this algorithm
  generates.
In the following, we let $\noff{A}$ denote the number of
  nonzero off-diagonal entries in the upper-triangular portion
  of $A$,
  and let $\dim{A}$ denote the dimension of $A$.

\begin{proposition}[Recursive Preconditioning]\label{pro:buildPrecon}
If $A_{0}$ is an irreducible \sddm-matrix, 
  and for each $i$ the 
   matrix $B_{i}$ is a $k$-ultra-sparsifier of $A_{i-1}$, then
\begin{enumerate}
\item [$(a)$]For $i \geq 1$, $\noff{A_{i}} \leq (3 h  /k) \noff{A_{i-1}}$.
\item [$(b)$]For $i \geq 1$, $\dim{A_{i}} \leq (2 h /k) \noff{A_{i-1}}$.

\item [$(c)$]For $i \geq 1$, $\dim{B_{i}} = \dim{A_{i-1}}$.

\item [$(d)$]Each of $B_{i}$ and $A_{i}$ is an irreducible  \sddm-matrix.
\item [$(e)$] $\ell \leq 2 \log_{4 h} n$.
\end{enumerate}
\end{proposition}
\begin{proof}
Let $n_{i}$ be the dimension of $A_{i}$.
Definition~\ref{def:ultraSparsifier} tells us that
\[
\noff{B_{i}} \leq n_{i-1} - 1 + h \noff{A_{i-1}} / k.
\]
Parts $(a)$, $(b)$, and $(d)$ now follow from Proposition~\ref{pro:chol}.
Part $(c)$ is obvious.  Part $(e)$ follows from part $(a)$.
\end{proof}  

Our recursive solver will use each matrix $B_{i}$ as a preconditioner
  for $A_{i-1}$.
But rather than solve systems in $B_{i}$ directly, it will reduce
  these to systems in $A_{i}$, which will in turn be solved recursively.
Our solver will use the preconditioned Chebyshev method, instead of the
  preconditioned conjugate gradient.
This choice is dictated by the requirements of our analysis
  rather than by common sense.
Our preconditioned Chebyshev method will not take the preconditioner
  $B_{i}$ as input.
Rather, it will take a subroutine $\Solve{B_{i}}$ that produces
  approximate solutions to systems in $B_{i}$.
So that we can guarantee that our solvers will be linear operators,
  we will fix the number of iterations that each will perform, as opposed to
  allowing them to terminate upon finding a sufficiently good solution.

For simplicity, we use the original Chebyshev iterative method~\cite{GolubVarga},
  as presented by Axelsson~\cite[Section 5.3]{Axelsson}.
While this variant is not numerically stable, it will not matter in this
  section in which we ignore numerical issues.
In particular, when one 
  assumes infinite precision this algorithm
  becomes identical to its stable variants.
In the next section, we will show that our use of the algorithm for a small
  number of iterations limits its instability.

\begin{figure}[ht]
\begin{algbox}
\noindent $\xx = \mathtt{precondCheby} (A, \bb, 
  f (\cdot ) , t, \lambda _{min}, \lambda _{max})$
\begin{tightlist}
\item [(0)] Set $\xx = \bvec{0}$ and $\rr = f (\bb)$.
\item [(1)] for $i = 1, \dotsc , t$,
\begin{tightlist}
\item [(a)] Set
$\theta_{i} = (2i-1) \pi / 2t$ and
$\tau_{i} = \left((\cos \theta_{i})(\lambda_{max} - \lambda_{min})/2 + (\lambda_{max} + \lambda_{min})/2\right)^{-1}$.
\item [(b)] Set $\xx =\xx - \tau_{i} f (A \xx) + \tau_{i} \rr $.
\end{tightlist}
\end{tightlist}
\end{algbox}
\end{figure}

\begin{proposition}[Linear Chebyshev]\label{pro:cheby}
Let $A$ be a positive definite matrix and $f$ be a positive definite, symmetric
  linear operator such that for some $\lambda_{max} \geq \lambda_{min} > 0$
\begin{equation}\label{eqn:pro:cheby}
  \lambda _{min} \inv{f} \pleq A \pleq \lambda _{max} \inv{f}.
\end{equation}
Let $\epsilon < 1$ and let
\begin{equation}\label{eqn:chebyParms}
  t \geq  \ceiling{\frac{1}{2} \sqrt{\frac{\lambda_{max}}{\lambda_{min}}} \ln \frac{2}{\epsilon }}.
\end{equation}
Then, the function 
  $\mathtt{precondCheby} (A, \bb  , f , t  , \lambda _{min}, \lambda _{max})$
 is a symmetric linear operator in $\bb$.
Moreover, if $Z$ is the matrix realizing this operator, then
\[
  (1-\epsilon ) \inv{Z} \pleq A \pleq (1+\epsilon ) \inv{Z}.
\]
\end{proposition}
\begin{proof}
An inspection of the pseudo-code reveals that the function computed
  by \texttt{precondCheby} can be expressed as a sum of monomials
  of the form $(f A)^{i} f$, from which it follows that this function
  is a symmetric linear operator.

Standard analyses of the preconditioned Chebyshev algorithm~\cite[Section 5.3]{Axelsson} imply
  that for all $\bb $,
\[
\norm{Z\bb - \inv{A} \bb }_{A} \leq \epsilon \norm{\inv{A}\bb }_{A}.
\]
Now, let $\lambda$ be any eigenvalue of $Z A$, let $\vv$ be the corresponding eigenvector,
  and let $\bb = A \vv$.
We then have
\[
  \epsilon \norm{\vv}_{A} \geq \norm{Z A \vv - \vv}_{A}
= 
\abs{\lambda -1} \norm{\vv}_{A}.
\]
So, $\abs{\lambda -1} \leq \epsilon$.
Applying Proposition~\ref{pro:pleq}, we 
  obtain
\[
  (1-\epsilon ) \inv{Z} \pleq A \pleq (1+\epsilon ) \inv{Z}.
\]
\end{proof}

We can now state the subroutine $\Solve{B_{i}}$ for $i = 1, \dotsc , \ell$.

\begin{algbox}
\noindent $\xx = \solve{B_{i}}{\bb}$
\begin{enumerate}

\item
 Set $\lambda_{min} = 1-2e^{-2}$, 
  $\lambda_{max} =  (1+2e^{-2}) k $ and $t = \ceiling{1.33\sqrt{k}}$,
  where $k$ is as set in \eqref{eqn:k} and in \texttt{BuildPreconditioners}
  for the system $A_{0}$.
\label{alg:solve:t}

\item  Set $\sss = L_{i}^{-1} P_{i}^{-1} \bb$. \label{alg:solvebi:s}
\item  Write $\sss = \left(\begin{array}{l}\sss_{0}\\
           \sss_{1} \end{array} \right)$,
  where the dimension of $\sss_{1}$ is the size of $A_{i}$.
  
\item Set $\yy_{0} = \sss_{0}$, and 

\begin{enumerate}
\item if $i = \ell $, set \label{alg:solvebi:ell}
 $\yy_{1} = Z_{\ell } \sss_{1}$
\item  else, set
  $\yy_{1} = \mathtt{precondCheby} (A_{i}, \sss_{1}, \Solve{B_{i+1}},
  t, \lambda_{min}, \lambda_{max}) $.  \label{alg:solvebi:cheby}
\end{enumerate}

\item Set $\xx =  P_{i}^{-T} L_{i}^{-T} 
  \left(\begin{array}{l}\yy_{0}\\  \yy_{1} \end{array} \right)$.
 \label{alg:solvebi:x}
\end{enumerate}
\end{algbox}

We have chosen the parameters $\lambda_{min}$, $\lambda_{max}$, and $t$ so that
  inequality~\eqref{eqn:chebyParms} holds 
  for $\epsilon = 2 e^{-2}$.
Our recursive algorithm only requires the solution of the systems $B_{i}$
  to some small constant error.
The constants given here are merely a simple choice that suffices.
It might be possible to obtain constant-factor improvements in running time
  by the choice of better constants.

We note that we apply $L_{i}^{-T}$ and $L_{i}^{-1}$ by forward 
  and backward substitution,
  rather than by constructing the inverses.

\begin{lemma}[Correctness of $\Solve{B_{i}}$]\label{lem:SolveBicorrect}
If $A$ is an irreducible  \sddm-matrix and
  $B_{i} \pleq A_{i-1} \pleq k B_{i}$ for all $i \geq 1$,
  then for $1 \leq i \leq \ell$,
\begin{dasenumerate}
\item The function $\Solve{B_{i}}$ is a symmetric linear operator.

\item The function 
  $\mathtt{precondCheby} (A_{i-1}, \bb  , \Solve{B_{i}}, t, \lambda_{min}, \lambda_{max}) $
  is a symmetric linear operator in $\bb$.

\item If $i \leq \ell - 1$ and $Z_{i}$ is the symmetric matrix such that
\[
Z_{i} \sss_{1} = 
   \mathtt{precondCheby} (A_{i}, \sss_{1} , \Solve{B_{i+1}}, t, \lambda_{min}, \lambda_{max}).
\]
Then,
\[
  (1-2e^{-2}) \inv{Z_{i}} \pleq A_{i} \pleq (1+2e^{-2}) \inv{Z_{i}}.
\]

\item 
\[
  (1-2e^{-2}) \inv{\Solve{B_{i}}}  \pleq B_{i} \pleq (1+2e^{-2})  \inv{\Solve{B_{i}}}.
\]

\end{dasenumerate}
\end{lemma}
\begin{proof}
We first prove $(a)$ and $(b)$ by 
  reverse induction on $i$.
The base case of our induction is when $i = \ell$, in which case
  \texttt{BuildPreconditioners} sets
  $Z_{\ell } = \inv{A_{\ell}}$, and so
\[
\Solve{B_{\ell}} = 
 P_{\ell}^{-T} L_{\ell}^{-T}
  \left(
     \begin{array}{ll}
       I & 0\\
       0 & Z_{\ell}
     \end{array}
   \right)
L_{\ell }^{-1} P_{\ell}^{-1} ,
\]
which is obviously a symmetric linear operator. 
Given that $\Solve{B_{i}}$ is a symmetric linear operator, 
  part $(b)$ for $A_{i-1}$ follows from
  Proposition~\ref{pro:cheby}.
Given that $(b)$ holds for $A_{i}$ and that
  the call to \texttt{precondCheby} is realized by a symmetric matrix $Z_{i}$, 
  we then have that
\[
\Solve{B_{i}} = 
 P_{i}^{-T} L_{i}^{-T}
  \left(
     \begin{array}{ll}
       I & 0\\
       0 & Z_{i}
     \end{array}
   \right)
L_{i}^{-1} P_{i}^{-1} 
\]
is a symmetric linear operator.
We may thereby establish that $(a)$ and 
  $(b)$ hold for all $\ell \geq i \geq 1 $.

We now prove properties $(c)$
  and $(d)$,
  again by reverse induction.
By construction $Z_{\ell} = \inv{A_{\ell }}$, so
  $(c)$ holds for $i = \ell$.
To see that if  $(c)$ holds for $i$,
  then $(d)$ does also, note that
\begin{align*}
  (1-2e^{-2}) \inv{Z_{i}} & \pleq A_{i} 
& \text{implies}\\
  (1-2e^{-2}) \inv{A_{i}} & \pleq Z_{i},
&  \text{by Proposition~\ref{pro:pleqPinv}, which implies}\\
  (1-2e^{-2})   
  \left(
     \begin{array}{ll}
       I & 0\\
       0 & \inv{A_{i}}
     \end{array}
   \right)
& \pleq 
  \left(
     \begin{array}{ll}
       I & 0\\
       0 &  Z_{i}
     \end{array}
   \right) 
&  \text{which implies}\\
\end{align*}
\begin{align*}
(1-2e^{-2}) \inv{B_{i}} & = 
(1-2e^{-2})   
 P_{i}^{-T} L_{i}^{-T}
  \left(
     \begin{array}{ll}
       I & 0\\
       0 & \inv{A_{i}}
     \end{array}
   \right) 
L_{i}^{-1} P_{i}^{-1} 
&  \text{(by Proposition~\ref{pro:chol} $(e)$)}
 \\
& \pleq  
 P_{i}^{-T} L_{i}^{-T}
  \left(
     \begin{array}{ll}
       I & 0\\
       0 &  Z_{i}
     \end{array}
   \right) 
L_{i}^{-1} P_{i}^{-1} \\
& = \Solve{B_{i}},
\end{align*}
which by Proposition~\ref{pro:pleqPinv} implies
 $(1-2e^{-2}) \inv{ \Solve{B_{i}}}
  \pleq 
  B_{i}
 $.
The inequality $B_{i} \pleq (1+2e^{-2}) \inv{\Solve{B_{i}}}$
  may be established similarly.

To show that when $(d)$ holds for $i$ 
  then $(c)$ holds for $i-1$,
  note that $(d)$ and
  $B_{i} \pleq A_{i-1} \pleq k \cdot B_{i}$
  imply
\[
  (1-2e^{-2}) \inv{\Solve{B_{i}}}
\pleq A_{i-1}
\pleq 
 k (1+2e^{-2}) \inv{\Solve{B_{i}}}.
\]
So, $(c)$ for $i-1$ now follows from Proposition~\ref{pro:cheby}
  and the fact that $\lambda_{min}, \lambda_{max}$ and $t$ have been
  chosen so that inequality \eqref{eqn:chebyParms} is satisfied with
  $\epsilon = 2 e^{-2}$.
\end{proof}

\begin{lemma}[Complexity of $\Solve{B_{i}}$]\label{lem:timeSolveBi}
If $A_{0}$ is a positive-definite
  irreducible,  
  $n$-by-$n$ \sddm-matrix with $2 m$ nonzero off-diagonal entries and 
  each $B_{i}$ is a $k$-ultra-sparsifier of $A_{i-1}$,
  then $\Solve{B_{1}}$ runs in time
\[
  O (n +  m).
\]
\end{lemma}
\begin{proof}
Let $T_{i}$ denote the running time of $\Solve{B_{i}}$.
We will prove by reverse induction on $i$ that
  there exists a constant $c$ such that
\begin{equation}\label{eqn:timeHyp}
  T_{i} \leq c \left(\dim{B_{i}} + (\gamma h + \delta)  (\noff{A_{i}} + \dim{A_{i}}) \right),
\end{equation}
where
\[
  \gamma = 196 \quad \text{and} \quad \delta = 15.
\]
This will prove the lemma as $\dim{B_{1}} = \dim{A_{0}} = n$,
  and Proposition~\ref{pro:buildPrecon} implies
\[
(\gamma  h + \delta ) (\noff{A_{1}} + \dim{A_{1}})
\leq 
(\gamma  h + \delta ) \frac{5 h m}{k}
\leq 
m \frac{5 \gamma  h^{2} + 5 \delta  h}{(14 h + 1)^{2}}
= O (m).
\]

To prove \eqref{eqn:timeHyp}, we note that there exists a constant $c$
  so that steps \ref{alg:solvebi:s}
  and \ref{alg:solvebi:x} take time at most $c (\dim{B_{i}})$
  (by Proposition~\ref{pro:chol}),
  step \ref{alg:solvebi:ell} takes time at most $c (\dim{A_{\ell}}^{2})$,
  and step \ref{alg:solvebi:cheby} takes time
  at most $t (c \cdot \dim{A_{i}}  + c \cdot \noff{A_{i}} + T_{i+1})$,
 where $t$ is as defined on step \ref{alg:solve:t} of $\Solve{B_{i}}$.

The base case of our induction will be $i = \ell$, in which
   case the preceding analysis implies 
\begin{align*}
T_{\ell} 
& \leq c \left(\dim{B_{\ell}} +  \dim{A_{\ell}}^{2}\right)
\\
& \leq c \left(\dim{B_{\ell}} +  (66 h + 6)\dim{A_{\ell}}\right),
& \text{ (by step \ref{alg:build:until} of \texttt{BuildPreconditioners})}
\end{align*}
which satisfies \eqref{eqn:timeHyp}.
We now prove \eqref{eqn:timeHyp} is true for $i< \ell$, assuming
  it is true for $i+1$.
We have
\begin{align*}
T_{i} 
& \leq 
  c \left(
     \dim{B_{i}}   \right)
    + t (c \cdot \dim{A_{i}} + c \cdot \noff{A_{i}} + T_{i+1})
\\
& \leq 
  c  \left[
     \dim{B_{i}}
    + t
  \big(
 \dim{A_{i}} + \noff{A_{i}} 
  + \dim{B_{i+1}} + (\gamma  h + \delta ) (\noff{A_{i+1}} + \dim{A_{i+1}})
 \big)
 \right]
\\
\intertext{(by the induction hypothesis)}
& \leq 
  c  \left[
     \dim{B_{i}}
    + t
  \big(
 2 \ \dim{A_{i}} + \noff{A_{i}} 
  + (\gamma  h + \delta ) (5 \ \noff{A_{i}} h / k )
 \big)
 \right]
\\
\intertext{(by Proposition~\ref{pro:buildPrecon})}
& \leq 
  c  \left[
     \dim{B_{i}}
    + t
  \left(
 2 \ \dim{A_{i}} + 6 \ \noff{A_{i}} 
 \right)
 \right],
\end{align*}
as $\gamma h^{2} + \delta h  \leq  k$.
As
\[
6 t \leq 6 \cdot (1.33 (14 h + 1) + 1)
\leq \gamma h + \delta ,
\]
we have proved that \eqref{eqn:timeHyp} is true for $i$ as well.
\end{proof}

We now state and analyze our ultimate solver.

\begin{algbox}
\noindent $\xx  = \solve{}{A, \bb, \epsilon }$
\begin{enumerate}
\item [1.]
Set $h = c_{3} \log_{2}^{c_{4}}\dim{A}$ (as in \eqref{eqn:h}) and $k = (14 h+1)^{2}$ (as in \eqref{eqn:k}).

Set $\lambda_{min} = 1-2e^{-2}$, 
  $\lambda_{max} =  (1+2e^{-2})k $ and $t = \ceiling{0.67 \sqrt{k} \ln (2/\epsilon )}$.

\item [2.]
  Run $\mathtt{BuildPreconditioners}(A)$.

\item [3.]
 $\xx = \mathtt{precondCheby} (A, \bb, \Solve{B_{1}},
  t, \lambda_{min}, \lambda_{max} )$.
\end{enumerate}
\end{algbox}

\begin{theorem}[Nearly Linear-Time Solver]\label{thm:solve}
On input an
  irreducible
  $n$-by-$n$ \sddm-matrix $A$ with $2 m$ nonzero off-diagonal entries and
  an $n$-vector  $\bb$,
  with probability at least $1 - 1/50$,
  $\Solve{} (A, \bb , \epsilon )$ runs in time
\[
O (m \log^{c_{4}} m  \log (1/\epsilon )) +
m \log^{c} m,
\]
where $c$ is some constant and $c_{4}$ is defined in \eqref{eqn:h},
  and produces an $\xxtil $ satisfying
\[
\norm{\xxtil - \inv{A} \bb }_{A} \leq \epsilon \norm{\inv{A} \bb }_{A}.
\]
\end{theorem}
\begin{proof}
By Proposition~\ref{pro:buildPrecon}, the numbers $\noff{A_{i}}$ are geometrically
  decreasing, and
  $\ell \leq 2 \log_{4 h} n$.
So we may use
  Theorem~\ref{thm:ultra} to show that
  the time required to build the preconditioners
  is at most $m \log^{O (1)} m$.
If each $B_{i}$ is a $k$-ultra-sparsifier of $A_{i-1}$,
  then the bound on the $A$-norm of the output
  follows by an analysis similar to that used to 
  prove  Lemma~\ref{lem:SolveBicorrect}.
In this case, we may use Lemma~\ref{lem:timeSolveBi} to
  bound on the running time of step $3$ 
  by
\[
  O\left(m t \right)
= O (m \sqrt{k} \log (1/\epsilon ))
 = O \left(m \log^{c_{4}} n \log (1/\epsilon ) \right).
\]
The probability that there is some $B_{i}$ that is not
  a $k$-ultra-sparsifier of $A_{i-1}$ is at most
\[
\sum_{i} \frac{1}{2 \ \dim{B_{i}}}
\leq 
\frac{
  \ell
}{ 2 (66 h +6)
}
\leq 
\frac{
  2 \log_{4 h} n
}{ 2 (66 h +6)
}
< 1/50,
\]
assuming $c_{3}, c_{4} \geq 1$.
\end{proof}

If the nonzero structure of $A$ is planar, then 
  by Theorem~\ref{thm:augmentTree},
  we can replace
  all the calls to \texttt{UltraSparsify} in the above algorithm with
  calls to \texttt{UltraSimple}.
By Theorem~\ref{thm:lowStretch}, this is like having
  $(k,h)$-ultra-sparsifiers with $h = O (\log n \log \log n)$.
Thus, the same analysis goes through with $h = O (\log n \log \log n)$,
  and the resulting linear system solver runs in time
\[
  O (n \log^{2} n + n \log n \  \log \log n \ \log (1/\epsilon )).
\]

\section{A Crude Stability Analysis}\label{sec:stability}

We will now show that the recursive solver described in Theorem~\ref{thm:solve} 
  works when all of the computations
  are carried out with limited precision.
In particular, we argue that $O (\log \kappa (A) \log^{c} n \log \epsilon^{-1})$
  bits of precision are sufficient.
While this bound is rather weak by the standards of Numerical Linear Algebra,
  it is sufficient for establishing many bounds on the asymptotic complexity of algorithms,
  such as those in~\cite{ChristianoEtAl,DaitchSpielman,kelnerMadryPropp,KelnerMillerPeng}.
We hope to one day see a better bound.
The main bottleneck in our analysis is that we have been unable to find
  a good analysis of the stability of the Preconditioned Chebyshev Method%
\footnote{While Golub and Overton~\cite{GolubOverton} suggest that such a stability 
  analysis should follow from the techniques they employ, the derivation of such a result
  is beyond the scope of the present paper.}.

As both the condition number and smallest eigenvalue of Laplacian matrices will play a substantial
  role in our analysis, we briefly relate these quantities to the weights of edges in the corresponding graph.
It follows from Ger{\u{s}}gorin's Circle Theorem that
  the largest eigenvalue of a \sdd-matrix is at most twice
  its largest diagonal entry
  (for the special case of Laplacians, see~\cite{AndersonMorley}).
A simple lower bound on the smallest eigenvalue of an irreducible SDDM matrix in terms of the
  lowest weight of an edge in the corresponding graph follows.

\begin{lemma}\label{lem:lminA}
Let $G$ be a connected weighted graph and let $A$ be either the Laplacian matrix of $G$
  or a principal square sub-matrix of the Laplacian.
Then the smallest nonzero eigenvalue of $A$ is at least $\min (8 w / n^{2}, w / n)$,
  where $w$ is the least weight of an edge of $G$ and $n$ is the dimension of $A$.
\end{lemma}
\begin{proof}
Fiedler~\cite{Fiedler} proved that the smallest eigenvalue of the Laplacian of 
  a connected, unweighted graph with $n$ vertices is at least
  $2 (1-\cos (\pi /n))$, which is at least $8 / n^{2}$ for $n \geq 2$.
In the weighted case, this implies that the smallest eigenvalue of the Laplacian
  is at least $8 w / n^{2}$ provided that all edge weights are at least $w$.

We now consider a sub-matrix of such a Laplacian.
Let $S$ be the set of vertices corresponding to the rows and columns of the sub-matrix.
The sub-matrix will have one diagonal block for each connected component of $S$.
Let $S_{1}$ be such a connected component.
The sub-matrix induced on $S_{1}$ can be decomposed into the sum of Laplacian matrix $A_{1}$ and
  a diagonal matrix $D_{1}$.
By the previous argument, the smallest nonzero eigenvalue of that Laplacian is at least
  $8 w / n^{2}$.
On the other hand, when we multiply the unit
  zero-eigenvector of that Laplacian by $D_{1}$
  we get $\bvec{1}^{T} D_{1} \bvec{1} / \sizeof{S_{1}}$.
The numerator equals the sum of the weights of edges on the boundary of $S_{1}$,
  which is at least $w$.
So, the smallest eigenvalue of the matrix induced on $S_{1}$ is at least
  $\min (8 w / n^{2}, w / n)$.
\end{proof}

We begin our analysis by asserting that the algorithm \texttt{UltraSparsify} is purely combinatorial
  and thus very stable.
It requires precision at most a polynomial in $n$ times the ratio of the
  largest to the smallest nonzero off-diagonal entry of its input.
This can be seen from an examination of the routines that it calls:
  \texttt{RootedUltraSparsify}, presented in Section~\ref{sec:ultra},
   and \texttt{Sparsify2} from~\cite{SpielmanTengSparsifier}.
In fact, the algorithm would barely suffer from rounding the weights of all edges
  in its input graph to powers of two.

We assume in the rest of this section that computations are performed with precision $u$,
 basing our analysis on those presented by Higham~\cite{Higham}.
To avoid the use of the  notation $O (u^{2})$, we employ Higham's~\cite{Higham} notation
\begin{equation}\label{eqn:gamma}
  \gamma_{j} = \frac{u j}{1 - u j}.
\end{equation}

We first address the issue that the matrices
  computed by \texttt{partialChol} will not be exactly the intended matrices by observing
  they are close enough to provide good preconditioners.

\begin{lemma}\label{lem:cholStable}
Let $L_{i}$ and $A_{i}$ be the matrices that would be output
  by $\mathtt{partialChol}$  on input $B_{i}$ if it were run with infinite precision,
  and let $\Lhat_{i}$ and $\Ahat_{i}$ be the matrices that are returned when it is run with precision $u$.
Let
\[
  \Bhat_{i}  = P \Lhat_{i} \left(  \begin{array}{ll}
 I & 0 \\
  0 & \Ahat _{i},
\end{array}  \right) \Lhat_{i}^{T} P^{T}.
\]
Then,
\[
  \left(1 - n \gamma_{n+1} \kappa (B_{i}) \right) B_{i}
\pleq
\Bhat_{i} 
\pleq 
   \left(1 - n \gamma_{n+1} \kappa (B_{i}) \right)^{-1} B_{i}.
\]
\end{lemma}
\begin{proof}[Proof Sketch]
Following the proof of  Lemma 2.1 in~\cite{DemmelCholesky}
(see also~\cite[Theorem 10.5]{Higham}),
we can show that every entry in $B_{i} - \Bhat_{i}$ is at most $\gamma_{n+1} \max_{j} B_{i} (j,j)$.
This implies that the norm of $B_{i} - \Bhat_{i}$ is at most $n \gamma_{n+1} \max_{j} B_{i} (j,j)$.
As $B_{i}$ is a positive semidefinite matrix, $\max_{j} B_{i} (j,j) \leq \lambda_{max} (B_{i})$.
So, for all $\xx$,
\[
\abs{  \frac{\xx^{T} \Bhat_{i} \xx - \xx^{T} B_{i} \xx}{\xx^{T} B_{i} \xx}}
\leq 
\frac{n \gamma_{n+1} \lambda_{max} (B_{i})}{\lambda_{min} (B_{i})}
= 
n \gamma_{n+1} \kappa (B_{i}).
\]
The lemma follows.
\end{proof}

We now prove that 
  the matrices produced by the routine \texttt{BuildPreconditioners} have condition
  numbers that are not too much larger than those of its input.

\begin{lemma}\label{lem:stabilityBP}
Let $A$ be a \sddm-matrix whose largest diagonal entry is between%
\footnote{The reason we make assumptions about the scale of $A$ is because we
  are bounding the condition number of a partial Cholesky factor.  
  We could avoid the need to make these assumptions if we instead computed partial $LDL^{T}$
  factorizations.
  The resulting algorithm would in fact be equivalent.}
 $1/2$ and $1$,
  and let $A_{1}, \dots , A_{\ell}$, $B_{1}, \dots , B_{\ell}$ and $L_{1}, \dots , L_{\ell}$
  be matrices produced by \texttt{BuildPreconditioners} 
  when it is run with precision $u$
  on input $A$.
If for each $i$ the matrix $B_{i}$ is a $k$-ultra-sparsifier of $A_{i-1}$, 
  and if
\[
  \gamma_{n+1} \leq \lambda_{min} (A) / 1000 n^{6},
\]
then
\begin{enumerate}
\item [a.] $ k^{l}  \leq n^{4}$, 

\item [b.] $\lambda_{min} (B_{i}) \geq \lambda_{min} (A) / 2 n^{4}$,
\item [c.] $\lambda_{max} (B_{i}), \lambda_{max} (A_{i})  \leq 3$,
\item [d.] $  \norm{L_{i}}_{\infty}  \leq 3 n, $ and
\item [e.] $ \norm{L_{i}^{-1}}_{\infty}  \leq 2 n^{3} / \sqrt{\lambda_{min} (A)}$.
\end{enumerate}
Finally,
  conditions (c) and (d) of Lemma~\ref{lem:SolveBicorrect}
  are satisfied by the matrices produced.
\end{lemma}
\begin{proof}
Proposition~\ref{pro:buildPrecon} tells us that
$\noff{A_{i}} \leq (3 h  /k) \noff{A_{i-1}}$.
As $\dim{A_{i}} \leq \noff{A_{i}}+1$ and \texttt{buildPreconditioners}
  stops when the dimension of $A_{i}$ goes below $66 h + 6$,
  we have
\[
  \ell \leq \log_{k/3h} 2m / (66h + 5) \leq \log_{k/3h} m.
\]
Since $k = (14h + 1)^{2}$, $k / 3 h \geq k^{1/2}$, which implies
\[
  k^{\ell} 
 \leq 
  k^{\log_{k^{1/2}} m}
  \leq m^{2}  
  \leq n^{4}.
\]

From the assumption that $B_{i}$ is a $k$-ultra-sparsifier of $A_{i-1}$,
  we know that
 $B_{i} \pleq A_{i-1} \pleq k B_{i}$, and so
\[
  \lambda_{max} (B_{i}) \leq \lambda_{max} (A_{i-1})
\]
and
\[
  \lambda_{min} (B_{i}) \geq (1/k) \lambda_{min} (A_{i-1}).
\]

Let
\[
  \Bhat_{i}  = P L_{i} \left(  \begin{array}{ll}
 I & 0 \\
  0 & A _{i},
\end{array}  \right) L_{i}^{T} P^{T},
\]
 and let $\epsilon$ be a number satisfying
\begin{equation}\label{eqn:stabilityBPeps}
  \epsilon \geq n \gamma_{n+1} \kappa (B_{i}).
\end{equation}
Lemma~\ref{lem:cholStable} then implies 
\begin{align*}
  \lambda_{max} (\Bhat _{i})  &
\leq (1-\epsilon)^{-1} \lambda_{max} (B_{i})
 \leq  (1 - \epsilon)^{-1} \lambda_{max} (A_{i-1}), \quad \text{and}
\\
  \lambda_{min} (\Bhat _{i}) &
  \geq (1-\epsilon) \lambda_{min} (B_{i})
 \geq (1 - \epsilon) (1/k) \lambda_{min} (A_{i-1}).
\end{align*}
As $A_{i}$ is a Schur complement of $\Bhat _{i}$, we have the 
  bounds~\cite[Corollary 2.3]{Zhang}
\[
  \lambda_{min} (\Bhat _{i}) \leq \lambda_{min} (A_{i}) \leq \lambda_{max} (A_{i}) \leq \lambda_{max} (\Bhat _{i}).
\]
Recall that $\lambda_{max} (A) \leq 2 \max_{i} A (i,i) \leq 2$.
If we set $\epsilon = 1 / 100 \ell$ and note that
  $(1-\epsilon)^{\ell} \geq e^{-.011}$, 
  we may now inductively prove that 
\[
  \lambda_{max} (B_{i}) \leq 2 (1-\epsilon)^{- (i-1)},
\quad 
  \lambda_{min} (B_{i}) \geq (1-\epsilon)^{i-1} \lambda_{min} (A) / k^{i},
\quad 
  \kappa (B_{i}) \leq 2 (1-\epsilon)^{-2 (i-1)} k^{i} / \lambda_{min} (A),
\]
and that \eqref{eqn:stabilityBPeps} is satisfied.
This establishes part $(b)$.
Part $(c)$ follows from
  $\lambda_{max} (A_{i}) \leq \lambda_{max} (\Bhat_{i}) \leq (1-\epsilon)^{-1} \lambda_{max} (B_{i})$.


Each matrix $L_{i}$ can be written in the form $R D$ where $R$ is a lower-triangular
  matrix with 1s on the diagonals and $D$ is a diagonal matrix of the form
\[
\left(
  \begin{array}{ll}
 D_{1} & 0 \\
   0 & I
\end{array}
 \right).
\]
The diagonal entries in $D_{1}$ are the square roots of the diagonal entries in $\Bhat _{i}$
  corresponding to nodes that are eliminated by \texttt{PartialChol}.
As $L_{i}$ is a Cholesky factor of a diagonally-dominant matrix, it is itself
  column diagonally-dominant, and thus $R$ is as well.
As every entry of $R$ is at most $1$, $\norm{R}_{\infty} \leq n$.
By a result of Malyshev~\cite{Malyshev} the same is true of $R^{-1}$.
As the diagonal entries of a symmetric matrix lie between its smallest and largest eigenvalues,
  the bounds we have proved on the eigenvalues of $\Bhat_{i}$ imply that its
  diagonals lie between $2 e^{3/100} $ and $ e^{-3/100} \lambda_{min} (A) / n^{4}$.
So,
\[
  \norm{L}_{\infty} \leq \norm{R}_{\infty} \norm{D}_{\infty} \leq 2 e^{3/100} n \leq 3 n,
\]
and
\[
  \norm{L^{-1}}_{\infty} \leq \norm{R^{-1}}_{\infty} \norm{D^{-1}}_{\infty} 
  \leq n \cdot n^{2} e^{-3/200} / \sqrt{\lambda_{min} (A)}
  \leq 2 n^{3} / \sqrt{\lambda_{min} (A)}.
\]

Finally, the proof that conditions (c) and (d) of Lemma~\ref{lem:SolveBicorrect} are satisfied
  depends upon inequality \eqref{eqn:chebyParms} being satisfied.
This inequality has a little bit of slack, and so it is satisfied even if
\[
  \lambda_{min} \geq (1- 2e^{-2}) / 1.003 \quad \text{and} \quad 
  \lambda_{max} \leq 1.003 (1 + 2e^{-2})k.
\]
Lemma~\ref{lem:cholStable} tells us that these conditions will hold if
\[
  n \gamma_{n+1} \kappa ( B_{i} ) \leq 0.003  ,
\]
and our assumptions on $\gamma_{n+1}$ guarantee that it does.
\end{proof}

The next lemma provides a very crude forward-error analysis of
  \texttt{precondCheby}.
We expect that it should be possible to obtain a tighter result for a stabler
  variant of the preconditioned Chebyshev method.

\begin{lemma}\label{lem:chebyStable}
Let $A$ be an $n$-by-$n$ \sddm \ matrix of norm at most $3$.
Let $B$ be a matrix such that $\lambda_{min} B \pleq A \pleq \lambda_{max}B$ with
  $\lambda_{min} \geq 1/2$ and $\lambda_{max} \leq 2 k$ for $k \geq 1$.
Let $\beta \geq 1$ be a number such that $\infnorm{B^{-1}} \leq \beta$. 
Assume there is a number $\theta$ and a procedure
 $f$ such that for all
  vectors $y$,
\[
  \infnorm{f (y) - B^{-1} y} \leq \theta \infnorm{y}.
\]
Also assume that $\theta$ and $u$ satisfy {\rm (}recall~\eqref{eqn:gamma}{\rm )}
\[
  1/10n \geq \theta \geq 12 k n \beta \gamma_{n}.
\]
Let $\xx^{t}$ be the result of running
  $\mathtt{precondCheby} (A, \bb, B^{-1}, t, \lambda_{min}, \lambda_{max} )$
  in infinite precision, and let
 $\xxhat^{t}$ be the result of running
  $\mathtt{precondCheby} (A, \bb, f, t, \lambda_{min}, \lambda_{max} )$
 with precision $u$.
Then
\[
  \infnorm{\xx^{t} - \xxhat^{t}} 
  \leq 
  17 \sqrt{n} (15 \beta  + 5k + 1)^{t}   \beta  \theta \infnorm{\bb}.
\]
\end{lemma}

\begin{proof}
Let $\xx^{i}$ be the vector computed in the $i$th iteration of
  \texttt{precondCheby} when it is run with infinite precision.
Also let $\xx^{0} = 0$ and let $b_{0} = \infnorm{\bb}$.
We have
\[
  \xx^{i+1} = \xx ^{i} - \tau_{i} B^{-1} (A \xx^{i}) + \tau_{i} B^{-1} \bb.
\]
Our conditions on  $\lambda_{min}$ and $\lambda_{max}$ imply
  that $\tau_{i} \leq 2$ for all $i$.
So,
\begin{align*}
  \norm{\xx^{i+1}} & \leq \norm{\xx^{i}} + 2 \norm{B^{-1} A} \norm{\xx^{i}} + 2 \norm{B^{-1} \bb}
\\
& \leq 5 k \norm{\xx^{i}} + 2 \norm{B^{-1} \bb}\\
& \leq (5 k + 1)^{i}  2 \norm{B^{-1} \bb}, & \text{by induction}.
\end{align*}

Using standard relations between the 2- and $\infty$-norms, we conclude
\[
  \infnorm{\xx^{i}} \leq 2 (5 k + 1)^{i-1} \sqrt{n} \beta b_{0}.
\]
In the rest of the proof, we make liberal use of relations such as these between norms,
  as well as the following inequalities:
\begin{align*}
  \infnorm{A} & \leq 6, \text{which follows from the diagonal dominance of $A$ and $\norm{A}\leq 3$, and}\\
\infnorm{B^{-1}} & \geq 1/6 \sqrt{n},
\text{which follows from $B \pleq 2A$ and $\infnorm{B^{-1}} \geq \norm{B^{-1}}/ \sqrt{n}$.}
\end{align*}

We now set
\begin{align*}
  \yy^{i} & = A \xx^{i}, \\
  \zz^{i} & = B^{-1} (\yy^{i}),
\end{align*}
and we let $\yyhat^{i}$ and $\zzhat^{i}$ be the analogous quantities computed
  using precision $u$ and the function $f$ instead of $B^{-1}$.

We compute
\begin{align*}
 \infnorm{ \yyhat^{i} - \yy^{i}}
& \leq 
  \gamma_{n} \infnorm{A} \infnorm{\xxhat^{i}}  
  + \infnorm{A} \infnorm{\xxhat^{i} - \xx^{i}} \quad  \text{(following \cite[Sec. 3.5]{Higham})}
\\
& \leq 
  \gamma_{n} \infnorm{A} \infnorm{\xx^{i}}  
  +   \gamma_{n} \infnorm{A} \infnorm{\xxhat^{i} - \xx^{i}}  
  + \infnorm{A} \infnorm{\xxhat^{i} - \xx^{i}}
\\
& \leq 
  6 \gamma_{n} \infnorm{\xx^{i}}  
  + 7 \infnorm{\xxhat^{i} - \xx^{i}},
\quad  \text{as $\infnorm{A} \leq 6$ and $\gamma_{n} \leq 1/6$.}
\end{align*}

We then compute
\begin{align*}
\infnorm{  \zzhat^{i} - \zz^{i}}
& =
\infnorm{  f (\yyhat^{i}) - B^{-1} \yy^{i}}
\\
& \leq 
\infnorm{  f (\yyhat^{i}) - B^{-1} \yyhat^{i}} + \infnorm{B^{-1} \yyhat^{i} - B^{-1} \yy^{i}}
& \text{by the triangle inequality}
\\
& \leq 
    \theta \infnorm{\yyhat^{i}}
+ \beta  \infnorm{\yyhat^{i} - \yy^{i}}
\\
& \leq 
    \theta \infnorm{\yy^{i}}
+ (\theta + \beta ) \infnorm{\yyhat^{i} - \yy^{i}}
\\
& \leq 
    \theta \infnorm{\yy^{i}}
+ 2 \beta  \infnorm{\yyhat^{i} - \yy^{i}}
& \text{as $\infnorm{B^{-1}} \geq 1/6 n \geq \theta  $}.
\end{align*}
If we now substitute our upper bound on $\infnorm{\yyhat^{i} - \yy^{i}}$
  and apply the inequality $\infnorm{\yy^{i}} \leq \infnorm{A} \infnorm{\xx^{i}}$,
  we obtain
\begin{align*}
\infnorm{  \zzhat^{i} - \zz^{i}}
& \leq 
\theta \infnorm{A} \infnorm{\xx^{i}}
+ 12 \beta   \gamma_{n} \infnorm{\xx^{i}}  
  + 14 \beta  \infnorm{\xxhat^{i} - \xx^{i}}
\\
& \leq 
 7 \theta  \infnorm{\xx^{i}}
  + 14 \beta  \infnorm{\xxhat^{i} - \xx^{i}},
\quad
\text{as $\theta \geq 12 \beta \gamma_{n}$ and $\infnorm{A} \leq 6$}.
\end{align*}
Finally, we find
\begin{align*}
  \infnorm{\xxhat^{i+1} - \xx^{i+1}}
& \leq 
(1 + \gamma_{3})
  \infnorm{\xxhat^{i} - \xx^{i}}
+
(1+\gamma_{3})  \infnorm{\zzhat^{i} - \zz^{i}}
+
\gamma_{3} \left( \infnorm{\zz^{i}} + \infnorm{\xx^{i}} + \beta b_{0}\right) 
+
\theta b_{0},
\end{align*}
where the terms involving $\gamma_{3}$ account for the imprecision introduced
 when computing the sum of the three vectors in the Chebyshev recurrence.
Using the upper bound on $\infnorm{\zzhat^{i} - \zz^{i}}$ and the inequality
  $\infnorm{\zz^{i}} \leq \infnorm{B^{-1} A} \infnorm{\xx^{i}} \leq 2
k \sqrt{n} \infnorm{\xx^{i}}$, we  
   see that this last expression
  is at most
\[
(1+\gamma_{3}) (1 + 14 \beta )   \infnorm{\xxhat^{i} - \xx^{i}}
+ 
\left((1 + \gamma_{3}) 7 \theta   + \gamma_{3} (1 + 2k\sqrt{n}) \right) \infnorm{\xx^{i}}
+ (\gamma_{3} \beta + \theta ) b_{0}.
\]
We simplify this expression by observing that
\[
(1+\gamma_{3}) (1 + 14 \beta )
\leq 
15 \beta ,
\]
\[
(1 + \gamma_{3}) 7 \theta   + \gamma_{3} (1 + 2k\sqrt{n}) 
\leq 
  8 \theta,
\]
and
\begin{align*}
  8 \theta \infnorm{\xx^{i}} +  (\gamma_{3} \beta + \theta ) b_{0}
& \leq 
16 \theta  (5 k + 1)^{i-1} \sqrt{n} \beta b_{0} + 2 \beta b_{0}
\quad \text{(as $\gamma_{3} \beta +\theta  \leq 2 \beta $)}
\\
& \leq 
17 \theta \sqrt{n} (5k + 1)^{i-1} \beta b_{0}.
\end{align*}
We conclude that
\[
  \infnorm{\xxhat^{i+1} - \xx^{i+1}}
\leq 
15 \beta 
  \infnorm{\xxhat^{i} - \xx^{i}}
+
17 \theta \sqrt{n} (5k + 1)^{i-1} \beta b_{0}.
\]
As $\xxhat^{0} = \xx^{0} = \bvec{0}$, we may apply this
  inequality inductively to obtain
\[
  \norm{  \xxhat^{t} - \xx^{t}}
\leq 
  17 (15 \beta  + 5k + 1)^{t} \theta \sqrt{n} \beta b_{0}.
\]

\end{proof}

We now establish analogous bounds on the stability of $\Solve{B_{i}}$.
\begin{lemma}\label{lem:stableSolveBi}
Assume that the largest diagonal entry of $A$ is between $1/2$ and $1$.
Let $A_{1}, \dots , A_{\ell}$, $B_{1}, \dots , B_{\ell}$ and $L_{1}, \dots , L_{\ell}$
  be as in Lemma~\ref{lem:stabilityBP}.
Let $i$ be less than $\ell$ and
 let $f_{i+1}$ be a function such that
  for all vectors $\bb$,
\[
  \infnorm{f_{i+1} (\bb) - \solve{B_{i+1}}{\bb}} \leq \theta \infnorm{\bb}.
\]
Let $\xx$ be the result of running $\Solve{B_{i}}$ on input $\bb$
  with full precision, 
and
  let $\xxhat$ be the result of running $\Solve{B_{i}}$ on input $\bb$
   with precision $u$,  using $f_{i+1}$ in place of $\Solve{B_{i+1}}$.
There exist constants $d_{1}, d_{2}, d_{3}$, and $d_{4}$ so that if
\[
  \gamma_{n} \leq \theta \lambda_{min}^{2} (A) / 50 n^{8}
\quad \text{and} \quad 
  \theta \leq 1/ d_{3} (n \kappa (A))^{d_{4}},
\]
then,
\[
    \infnorm{\xxhat - \xx}  \leq 
  \theta  d_{1} (n \kappa (A))^{d_{2} \sqrt{k}} \infnorm{\xx}.
\]
The same bound holds for $i = \ell$ if $\Solve{B_{\ell}}$ is executed
  using precision $u$, assuming that $Z_{\ell}$ is applied by forward
  and backward substitution through Cholesky factors of $A_{\ell}$.
\end{lemma}
\begin{proof}
We first consider the case when $i < \ell$.
Let $\ss, \ss_{0}, \ss_{1}$ and 
  $\yy$ denote the vectors computed by  $\Solve{B_{i}}$
  when it runs in full precision.
Similarly, let $\sshat, \sshat_{0}, \sshat_{1}$ and $\yyhat$ be the corresponding vectors
  computed when $\Solve{B_{i}}$ runs with precision $u$
  and uses $f$.
By inequality (8.2) of Higham~\cite{Higham}, we have
\begin{equation}\label{eqn:stableSolveBiCond}
  \infnorm{\sshat - \ss} \leq 
  \frac{\mathrm{Cond} (L_{i})\gamma_{n}}
       {1-\mathrm{Cond} (L_{i})\gamma_{n}} \infnorm{\ss},
\end{equation}
where 
\[
\mathrm{Cond} (L_{i}) \leq \infnorm{L_{i}} \infnorm{L_{i}^{-1}}.
\]
By Lemma~\ref{lem:stabilityBP}, we know that this product is
   at most $6 n^{4} / \sqrt{\lambda_{min} (A)}$.
We have assumed $\gamma_{n}$ is small enough so that this gives
\[
  \infnorm{\sshat - \ss} \leq 
7 n^{4} \gamma_{n} \infnorm{\ss}  / \sqrt{\lambda_{min} (A)}.
\]
We also know that
\[
  \infnorm{\ss} \leq \infnorm{L_{i}^{-1}} \infnorm{\bb} 
  \leq 2 n^{3} \infnorm{\bb}  / \sqrt{\lambda_{min} (A)},
\]
and so
\[
  \infnorm{\sshat_{1} - \ss_{1}} \leq 
  \infnorm{\sshat - \ss} \leq 
  12 n^{7} \gamma_{n} \infnorm{\bb}  / {\lambda_{min} (A)}.
\]
Again applying our assumptions on $\gamma_{n}$, we derive
\[
\infnorm{\sshat_{1}}
 \leq   \infnorm{\sshat} \leq 3 n^{3} \infnorm{\bb}  / \sqrt{\lambda_{min} (A)}.
\]

We now examine the relationship between $\yy$ and $\yyhat$.
For $i < \ell$, Lemma~\ref{lem:chebyStable} tells us that
\begin{align*}
  \infnorm{\yyhat - \yy} 
  & \leq \alpha \theta \infnorm{\sshat_{1}} + \infnorm{\Solve{B_{i+1}}} \infnorm{\sshat_{1} - \ss_{1}},
\end{align*}
where
\[
  \alpha =   17 (15 \beta  + 5k + 1)^{t} \theta \sqrt{n} \beta,
\]
and
\[
  \beta = \infnorm{\Solve{B_{i+1}}} \leq \sqrt{n} \norm{\Solve{B_{i+1}}}.
\]
By part $(d)$ of Lemma~\ref{lem:SolveBicorrect} and by Lemma~\ref{lem:stabilityBP}, we know
\[
  \norm{\Solve{B_{i+1}}} \leq 2 (1+2 e^{-2}) \norm{B_{i+1}^{-1}} \leq 
  6 n^{4} / \lambda_{min} (A).
\]
From the specification of $\Solve{B_{i}}$, we have $t = \ceiling{1.33 \sqrt{k}}$.
So,
\[
  \alpha \leq d_{1} (n \kappa (A))^{d_{2} \sqrt{k}},
\]
for some constants $d_{1}$ and $d_{2}$.
From our assumption on the relation between $\theta$ and $\gamma_{n}$ we conclude
\[
  \infnorm{\yyhat - \yy} 
 \leq 2 \alpha \theta \infnorm{\sshat}.
\]
We also know that
\[
  \infnorm{\yy} \leq \infnorm{\Solve{B_{i+1}}} \infnorm{\ss}
 \leq 
    12 n^{7.5}  \infnorm{\bb} / \lambda_{min}^{3/2} (A) ,
\]
and by our assumptions on $\theta$ that
\[
  \infnorm{\yyhat}
 \leq 
    13 n^{7.5}  \infnorm{\bb} / \lambda_{min}^{3/2} (A).
\]

There are two sources of discrepancy between $\xxhat$ and $\xx$:
  error from multiplying $\yyhat $ by $L_{i}^{-T}$ and error from
  the difference between $\yyhat$ and $\yy$.
Let $\xxtil = L_{i}^{-T} \yyhat$.
We have
\[
  \infnorm{\xxhat - \xx} \leq \infnorm{\xxhat -\xxtil} + \infnorm{\xxtil - \xx}.
\]

As with the derivation following \eqref{eqn:stableSolveBiCond},
  we find
\begin{align*}
  \infnorm{\xxhat - \xxtil} 
& \leq 7 n^{4} \gamma_{n} \infnorm{\xxtil}  / \sqrt{\lambda_{min} (A)}
\\
& \leq 7 n^{4} \gamma_{n} \infnorm{L_{i}^{-T}} \infnorm{\yyhat }  / \sqrt{\lambda_{min} (A)}.
\end{align*}
We also have
\begin{align*}
\infnorm{\xxtil - \xx}
\leq 
\infnorm{L_{i}^{-T}} \infnorm{\yyhat - \yy}.
\end{align*}
We conclude that
\[
    \infnorm{\xxhat - \xx}  \leq 
  \theta  d_{1} (n \kappa (A))^{d_{2} \sqrt{k}} \infnorm{\bb},
\]
for some other constants $d_{1}$ and $d_{2}$.
Finally, we prove the main claim of the lemma for some different constants $d_{1}$ and $d_{2}$
  by observing that
\[
\infnorm{\bb} \leq \infnorm{\Solve{B_{i}}^{-1}}  \infnorm{\xx},
\]
and that
\[
  \infnorm{\Solve{B_{i}}^{-1}} \leq
  \sqrt{n} \norm{\Solve{B_{i}}^{-1}} \leq 
 (1-2e^{-2})^{-1} \sqrt{n} \norm{B_{i}}
 \leq 5 \sqrt{n},
\]
where the second inequality follows from part $(d)$ of Lemma~\ref{lem:SolveBicorrect}
  and the last inequality follows from part $(c)$ of Lemma~\ref{lem:stabilityBP}.

In the case $i = \ell$, we may apply the same analysis but 
  without the need to account for a call to \texttt{precondCheby}.
While there is no function $f_{i+1}$ (or we can view it as the identity),
  we still require the upper bound induced on $\gamma_{n}$ through $\theta$.
\end{proof}

\begin{theorem}\label{thm:stable}
Let $A$ be an \sddm -matrix whose largest diagonal entry is between $1/2$ and $1$.
Let $\epsilon < 1/2$, 
  let $\xx$ be the result of running   $\Solve{} (A, \bb , \epsilon )$ 
  in infinite precision,
  and let  $\xxhat$ be the result of running this routine and all of its subroutines
  with precision $u$.
There exist constants $d_{1}, d_{2}$ and $d_{3}$ so that if
\[
  u \leq (\epsilon / d_{1} \kappa (A) )^{d_{2} \log^{d_{3}} n},
\]
then
\[
\infnorm{\xx - \xxhat } \leq \epsilon \infnorm{\xx}.
\]
\end{theorem}
That is, it suffices to run $\Solve{}$ with $O (\log (\kappa (A)/ \epsilon ) \log^{d_{3}} n)$
  bits of precision.

\begin{proof}
As  $\Solve{}$ uses \texttt{precondCheby} to solve systems in $A$
  by using $B_{1}$ as a preconditioner,
  we begin by examining the requirements on $\Solve{B_{1}}$.

For each $i$, let $f_{i} (\bb_{i})$ be the output of $\Solve{B_{i}}$ when it and all of  
  its subroutines are executed with precision $u$ on input $\bb_{i}$.
We need to establish that there exist numbers $\theta_{i}$ that satisfy
  the conditions of Lemma~\ref{lem:stableSolveBi} for $f_{i}$
  (where in the statement of that lemma $\Solve{B_{i}}$ is treated as the infinite-precision
  linear operator).
We may assume that the constants $d_{1}, d_{2}, d_{3}$ and $d_{4}$ from Lemma~\ref{lem:stableSolveBi}
  are all at least $1$.
So, it suffices to choose $\theta_{1}$ less than $1/d_{3} (n \kappa (A))^{d_{4}}$,
\[
  \theta_{i} = \frac{\theta_{1}}{d_{1} (n \kappa (A))^{d_{2} \sqrt{k} (i-1)}},
\]
and $u$ so that
\[
  \gamma_{n} \leq \theta_{\ell} \lambda_{max}^{2} (A) / 50 n^{8}.
\]

If we now apply the analysis of \texttt{precondCheby} from Lemma~\ref{lem:chebyStable},
  assuming that the constants $d_{3}$ and $d_{4}$ are big enough that
  $1/d_{3} (n \kappa (A))^{d_{4}} < 1/10n$,
  we see that
\[
\infnorm{\xx - \xxhat}
\leq 
  17 \sqrt{n} (15 \beta  + 5k + 1)^{t}   \beta  \theta_{1} \infnorm{\bb},
\]
where
\[
  \beta = \infnorm{\Solve{B_{1}}},
\]
and
\[
t \leq \sqrt{k} \ln (2 / \epsilon).
\]
As in the proof of Lemma~\ref{lem:stableSolveBi}, we can show that
\[
\infnorm{\Solve{B_{1}}} \leq 
    6 n^{4} / \lambda_{min} (A)
\]
and
\[
  \infnorm{\bb}   \leq 5 \sqrt{n} \infnorm{\xx}.
\]
To finish the proof, we observe that $\ell \leq \log_{2} n$
  and recall from \eqref{eqn:k} that
  $\sqrt{k} = (14 h + 1)$ and that $h$ was set in \eqref{eqn:h} to
  $c_{3} \log_{2}^{c_{4}} n$, where $c_{3}$ and $c_{4}$
  were parameters related to the quality of the ultra-sparsifiers that
  we could produce.
We thus conclude that we can find (different) constants $d_{1}, d_{2}$ and $d_{3}$
  that satisfy the claims of the theorem.
\end{proof}

\section{Computing Approximate Fiedler Vectors}\label{sec:eigs}

Fiedler~\cite{Fiedler} was the first to recognize that 
  the eigenvector associated with  the second-smallest eigenvalue
  of the Laplacian matrix of a graph could be used to partition
  a graph.  
From a result of Mihail~\cite{Mihail89}, we know that any
  vector whose Rayleigh quotient is close to this eigenvalue
  can also be used to find a good partition.
We call such a vector an \textit{approximate Fiedler vector}.

\begin{definition}[Approximate Fiedler Vector]\label{def:fiedler}
For a Laplacian matrix $A$, $\vv$ is an $\epsilon$-approximate Fiedler vector
  if $\vv$ is orthogonal to the all-1's vector and
\[
  \frac{\vv^{T} A \vv}{\vv^{T} \vv } \leq (1+\epsilon) \lambda_{2} (A),
\]
where $\lambda_{2} (A)$ is the second-smallest eigenvalue of $A$.
\end{definition}

Our linear system solvers may be used to 
  quickly compute $\epsilon$-approximate Fiedler vectors.
The algorithm \texttt{ApproxFiedler} does so
 with probability at least $1-p$. 
This algorithm works by applying the inverse power method to a 
  random initial vector and by using $\Solve{}$ to accomplish the inversion.
As there is some chance that this algorithm could fail, or that \texttt{BuildPreconditioners}
  could fail, we run the entire process $\log_{2} 1/p$ times and return the best result.

\begin{algbox}
\noindent $\vv  = \texttt{ApproxFiedler}(A, \epsilon , p)$\\
\begin{enumerate}
\item Set $\lambda_{min} = 1-2e^{-2}$, 
  $\lambda_{max} =  (1+2e^{-2})k $ and $t = \ceiling{0.67 \sqrt{k} \ln (8 /\epsilon )}$.

\item Set $s = 8 \ln (18 (n-1)/ \epsilon ) / \epsilon $.

\item For $a = 1, \dots , \ceiling{\log_{2} 1/p}$.
\begin{dasenumerate}

\item 
  Run $\mathtt{BuildPreconditioners}(A)$.

\item Choose $\rr^{0} $ to be a uniform random unit vector orthogonal to the all-1's vector.

\item For $b = 1, \dots , s$
\begin{itemize}
\item [] $\rr^{b}  =  \mathtt{precondCheby} (A, \rr^{b-1} , \Solve{B_{1}},
  t, \lambda_{min}, \lambda_{max} )$.
\end{itemize}
\item Set $\vv_{a} = \rr^{s}$.
\end{dasenumerate}
\item Let $a_{0}$ be the index of the vector minimizing
  $\vv_{a_{0}}^{T} A \vv_{a_{0}} / \vv_{a_{0}}^{T} \vv_{a_{0}}$.
\item Set $\vv  = \vv_{a_{0}}$.
\end{enumerate}
\end{algbox}

\begin{theorem}\label{thm:approxFiedler}
On input a Laplacian matrix $A$ 
  with $m$ nonzero entries and $\epsilon , p > 0$, with probability at least $1-p$,
  $\mathtt{ApproxFiedler} (A, \epsilon , p )$  computes an
  $\epsilon$-approximate Fiedler vector of $A$ in time
\[
 O ( m \log^{c} m \log (1/p) \log (1/\epsilon ) / \epsilon ),
\]
for some constant $c$.
\end{theorem}

Our proof of Theorem~\ref{thm:approxFiedler} will use the following
  proposition.

\begin{proposition}\label{pro:approxFiedler}
If $Z$ is a matrix such that 
\[
  (1-\epsilon) \pinv{Z} \pleq A \pleq (1+\epsilon) \pinv{Z},
\]
and $\vv$ is a vector orthogonal to the all-1's vector such that
  $\vv^{T} \pinv{Z} \vv \leq (1+\epsilon) \lambda_{2} (\pinv{Z})$,
  for some $\epsilon \leq 1/5$,
then
  $\vv$ is a $4\epsilon$-approximate Fiedler vector of $A$.
\end{proposition}
\begin{proof}
We first observe that
\[
  \lambda_{2} (\pinv{Z}) \leq \lambda_{2} (A) /   (1-\epsilon).
\]
We then compute
\begin{align*}
\vv^{T} A \vv 
& \leq (1+\epsilon ) \vv^{T} \pinv{Z} \vv \\
& \leq (1+\epsilon ) (1+\epsilon  ) \lambda_{2} (\pinv{Z})\\
& \leq (1+\epsilon ) (1+\epsilon  ) \lambda_{2} (A) /   (1-\epsilon)\\
& \leq 
 (1+4 \epsilon) \lambda_{2} (A),
\end{align*}
for $\epsilon \leq 1/5$.
\end{proof}

\begin{proof}[Proof of Theorem~\ref{thm:approxFiedler}]
As we did in the proof of Lemma~\ref{lem:SolveBicorrect} and 
  Theorem~\ref{thm:solve},
  we can show that 
  $\mathtt{precondCheby} (A, \bb , \Solve{B_{1}},
  t, \lambda_{min}, \lambda_{max} )$
  is a linear operator in $\bb$.
Let $Z$ denote the matrix realizing this operator.
As in the proof of Lemma~\ref{lem:SolveBicorrect}, we can show that 
  $(1-\epsilon / 4) \pinv{Z} \pleq A \pleq (1+\epsilon/ 4) \pinv{Z}$.

By Proposition~\ref{pro:approxFiedler},
  it suffices to show that with probability at least $1/2$
 each vector $\vv_{a}$ satisfies
\[
  \vv_{a}^{T} \pinv{Z} \vv_{a} / \vv_{a}^{T} \vv_{a} \leq (1 + \epsilon /4) \lambda_{2} (\pinv{Z}).
\]
To this end,
  let $0 = \mu_{1} \leq \mu_{2} \leq \dotsb \leq \mu_{n}$
  be  the eigenvalues of $\pinv{Z}$, and let
  $\bvec{1} = \uu_{1}, \dotsc , \uu_{n}$ be corresponding eigenvectors.
Let
\[
  \rr^{0}  = \sum_{i \geq 2} \alpha_{i} \uu_{i},
\]
and recall that (see \textit{e.g.} \cite[Lemma B.1]{SankarSpielmanTeng})
\[
\prob{}{\abs{\alpha_{2}} \geq  2 / 3 \sqrt{(n-1)}} \geq 
 \frac{2}{\sqrt{2 \pi}} \int_{2/3}^{\infty} e^{-t^{2}/2} \ dt \geq 0.504.
\]
Thus, with probably at least $1/2$, the call to \texttt{BuildPreconditioners}
  succeeds and $\abs{\alpha_{2}} \geq  2 / 3 \sqrt{(n-1)}$.
In this case, 
\begin{equation}\label{eqn:fiedlerk}
  s \geq  8 \ln (8/ \alpha_{2}^{2} \epsilon ) / \epsilon.
\end{equation}
We now show that this inequality
  implies that 
  $\rr^{s}$ satisfies
\[
  \frac{(\rr^{s})^{T} \pinv{Z} \rr^{s}}{(\rr^{s})^{T} \rr^{s}}
\leq 
  (1 + \epsilon / 4) \mu_{2}.
\]
To see this, let $j$ be the greatest index such that
  $\mu_{j} \leq (1+\epsilon / 8) \mu_{2}$, and
 compute
\[
  \rr^{s} = Z^{s} \rr^{0} = \sum_{i \geq 2} \uu_{i} \alpha_{i} / \mu_{i}^{s} ,
\]
so
\begin{align*}
 \frac{(\rr^{s})^{T} \pinv{Z} \rr^{s}}{(\rr^{s})^{T} \rr^{s}}
& = 
\frac{
\sum_{i \geq 2} \alpha_{i}^{2} / \mu_{i}^{2s - 1}
}{
\sum_{i \geq 2} \alpha_{i}^{2} / \mu_{i}^{2s}
}\\
& \leq 
\frac{
\sum_{j \geq i \geq 2} \alpha_{i}^{2} / \mu_{i}^{2s - 1}
}{
\sum_{j \geq i \geq 2} \alpha_{i}^{2} / \mu_{i}^{2s}
}
+
\frac{
\sum_{i > j} \alpha_{i}^{2} / \mu_{i}^{2s - 1}
}{
\sum_{i \geq 2} \alpha_{i}^{2} / \mu_{i}^{2s}
}\\
& \leq 
\mu_{j}
+
\frac{
\sum_{i > j} \alpha_{i}^{2} / \mu_{i}^{2s - 1}
}{
\alpha_{2}^{2} / \mu_{2}^{2s}
}\\
& \leq 
(1+\epsilon / 8) \mu_{2}
+
\mu_{2}
\left(
\frac{
\sum_{i > j} \alpha_{i}^{2} (\mu_{2}/\mu_{i})^{2s - 1}
}{
 \alpha_{2}^{2}
}
 \right)
\\
& \leq 
(1+\epsilon/8) \mu_{2}
+
\mu_{2}
\left(
\frac{
\sum_{i > j} \alpha_{i}^{2} (1/ (1+\epsilon/8 ))^{2s - 1}
}{
 \alpha_{2}^{2}
}
 \right)
\\
& \leq 
(1+\epsilon / 8) \mu_{2}
+
\mu_{2}
\sum_{i > j} \alpha_{i}^{2} \epsilon /8
& \text{(by inequality \eqref{eqn:fiedlerk})}
\\
& \leq 
(1+\epsilon/8) \mu_{2}
+
\mu_{2}
(\epsilon /8)
\\
& \leq 
(1 + \epsilon / 4) \mu_{2}.
\end{align*}
\end{proof}

\section{Laplacians and Weighted Graphs}\label{sec:graphs}
We will find it convenient to describe and analyze our preconditioners
  for Laplacian matrices in terms of weighted graphs.
This is possible because of the isomorphism between Laplacian matrices
  and weighted graphs.
To an $n$-by-$n$ Laplacian matrix $A$, we associate the graph with
  vertex set $\setof{1,\dotsc ,n}$ having an edge between vertices 
  $u$ and $v$ of weight $-A (u,v)$ for each $u$ and $v$ such that
  $A (u,v)$ is nonzero.

All the graphs we consider in this paper will be weighted.
If $u$ and $v$ are distinct vertices in a graph, we write
  $\edg{u,v}$ to denote an edge between $u$ and $v$
  of weight 1.
Similarly, if $w > 0$, then we write
  $w \edg{u,v}$ to denote an edge between $u$
  and $v$ of weight $w$.
A weighted graph is then a pair $G = (V,E)$ where $V$ is a set of
  vertices and $E$ is a set of weighted edges on $V$, each of which
  spans a distinct pair of vertices.
The Laplacian matrix $L_{G}$ of the graph $G$ is the matrix
  such that
\[
  L_{G} (u,v) = 
\begin{cases}
-w  & \text{if there is an edge $w\edg{u,v} \in E$}
\\
0  & \text{if $u \not = v$ and there is no edge between $u$ and $v$ in $E$}\\
\sum_{w\edg{u,x} \in E} w & \text{if $u = v$}.
\end{cases}
\]

We recall that for every vector $\xx \in \Reals{n}$,
\[
  \xx^{T} L_{G} \xx = \sum_{w \edg{u,v} \in E} w (\xx_{u} - \xx_{v})^{2}.
\]

For graphs $G$ and $H$ on the same set of vertices, we define the graph $G + H$ to be the graph
  whose Laplacian matrix is $L_{G} + L_{H}$.

\section{Graphic Inequalities, Resistance, and Low-Stretch Spanning Trees}\label{sec:gi}

In this section, we introduce the machinery of ``graphic inequalities''
  that underlies the proofs in the rest of the paper.
We then introduce low-stretch spanning trees, and use graphic
  inequalities to bound how well a low-stretch spanning tree
  preconditions a graph.
This proof provides the motivation for the construction
  in the next section.

We begin by overloading the notation $\pleq$
  by writing
\[
  G \pleq H  \qquad \text{or} \qquad E \pleq F
\]
if $G = (V,E)$ and $H = (V,F)$ are two graphs such that
  their Laplacian matrices, $L_{G}$ and $L_{H}$ satisfy
\[
L_{G} \pleq L_{H}.
\]

Many facts that have been used in the chain of work related to
 this paper can be simply expressed with this notation.
For example, the Splitting Lemma of \cite{SupportGraph}
  becomes
\[
  A_{1} \pleq B_{1} \quad 
\text{and} \quad 
  A_{2} \pleq B_{2} \quad 
\text{implies} \quad 
A_{1} + A_{2} \pleq B_{1} + B_{2}.
\]
We also observe that if $B$ is a subgraph of $A$,
  then
\[
B \pleq A.
\]

We define the \textit{resistance}
  of an edge to be the reciprocal of its weight.
Similarly, we define the resistance of a simple path to
  be the sum of the resistances of its edges.
For example, the resistance of the path
  $w_{1}\edg{1,2}, w_{2}\edg{2,3}, w_{3}\edg{3,4}$
  is $(1/w_{1} + 1/w_{2} + 1/w_{3})$.
Of course, the resistance of a trivial path with one vertex and no edges
  is zero.
If one multiplies all the weights of the edges in a path by $\alpha$,
  its resistance decreases by a factor of $\alpha$.

The next lemma says that a path of resistance $r$
  supports an edge of resistance $r$.
This lemma may be derived from the 
  Rank-One Support Lemma of~\cite{SupportTheory},
  and appears in simpler form as the
  Congestion-Dilation Lemma of~\cite{SupportGraph}
  and Lemma 4.6 of~\cite{Gremban}.
We present a particularly simple proof.

\begin{lemma}[Path Inequality]\label{lem:pathResistance}
Let $e = w \edg{u,v}$
  and let $P$ be a path from $u$ to $v$.
Then,
\[
e \pleq w \ \res{P} \cdot P.
\]
\end{lemma}
\begin{proof}
After dividing both sides by $w$, it suffices to consider the case $w = 1$.
Without loss of generality, we may assume that $e = \edg{1,k+1}$
   and that
  $P$ consists of the edges $w_{i}\edg{i,i+1}$ for $1 \leq i \leq k$.
In this notation, the lemma is equivalent to
\[
  \edg{1,k+1}
 \pleq 
\left(\sum_{i} \frac{1}{w_{i}} \right)
\big(  w_{1}\edg{1,2} +  w_{2}\edg{2,3} + \dots + w_{k}\edg{k,k+1} \big).
\]
We prove this for the case $k = 2$.
The general case follows by induction.

Recall Cauchy's inequality, which says that for all $0 < \alpha < 1$,
\[
(a+b)^{2} \leq a^{2} / \alpha + b^{2} / (1-\alpha ).
\]
For $k = 2$, the lemma is equivalent to
 the statement that for all $\xx \in \Reals{3}$,
\[
 \xx^{T} L_{\edg{1,3}} \xx \leq \left(\frac{1}{w_{1}} + \frac{1}{w_{2}} \right)
  \xx^{T} \left( L_{w_{1} \edg{1,2}} 
  +  L_{w_{2} \edg{2,3}}  \right) \xx.
\]
This is equivalent to
\[
  (x_{1} - x_{3})^{2} \leq (1 + w_{1}/w_{2}) (x_{1} - x_{2})^{2}
                         + (1 + w_{2}/w_{1}) (x_{2} - x_{3})^{2},
\]
which follows from Cauchy's inequality with $\alpha = w_{2}/(w_{1} + w_{2})$.
\end{proof}

Recall that a \textit{spanning tree} of a weighted graph $G = (V,E)$
  is a connected subgraph of $G$ with exactly $\sizeof{V}-1$ edges.
The weights of edges that appear in a spanning tree are assumed to be
  the same as in $G$.
If $T$ is a spanning tree of a graph $G = (V,E)$, then
  for every pair of vertices $u,v \in  V$,
  $T$ contains a unique path from $u$ to $v$.
We let $T (u,v)$ denote this path.
We now use graphic inequalities to 
  derive a bound on how well $T$ preconditions $G$.
This bound strengthens a bound 
  of Boman and Hendrickson~\cite[Lemma~4.9]{SupportTheory}.

\begin{lemma}[Tree Preconditioners]\label{lem:treePrecon}
Let $G = (V,E)$ be a graph and let $T$ be a spanning tree of $G$.
Then,
\[
T \pleq G \pleq 
\left(\sum_{e \in E} \frac{\res{T (e)}}{\res{e}} \right) \cdot T.
\]
\end{lemma}
\begin{proof}
As $T$ is a subgraph of $G$, $T \pleq G$ is immediate.
To prove the right-hand inequality, we compute
\begin{align*}
E & =
\sum_{e \in E} e
\\
& \pleq 
\sum_{e \in E} 
\frac{\res{T (e)}}{\res{e}}\cdot 
T (e),
& \text{by Lemma~\ref{lem:pathResistance}}
\\
& \pleq 
\left(\sum_{e \in E} 
\frac{\res{T (e)}}{\res{e}} \right)
\cdot T,
&
\text{as $T (e) \pleq T$.}
\end{align*}
\end{proof}

\begin{definition}[Stretch]\label{def:stretch}
Given a tree $T$ spanning a set of vertices $V$ and a weighted
  edge $e = w \edg{u,v}$ with $u,v \in V$,
  we define
  the \textit{stretch} of $e$ with respect to $T$ to be
\[
  \stretch{T}{e} = \frac{\res{T (e)}}{\res{e}}
=
 w \cdot \res{T (e)}
.
\]
If $E$ is a set of edges on $V$, then we define
\[
  \stretch{T}{E} = \sum_{e \in E}\stretch{T}{e}.
\]
\end{definition}

With this definition, the statement of Lemma~\ref{lem:treePrecon} may be simplified
  to
\begin{equation}\label{eqn:treeStretch}
T \pleq   G \pleq \stretch{T}{E} \cdot T.
\end{equation}

We will often use the following related inequality, which follows
  immediately from Lemma~\ref{lem:pathResistance} and the definition of stretch.
\begin{equation}\label{eqn:edgeStretch}
w \edg{u,v} \pleq \stretch{T}{w \edg{u,v}} \ T (u,v) = w \ \stretch{T}{\edg{u,v}} \  T (u,v).
\end{equation}

\section{Preconditioning with Augmented Low-Stretch Trees}\label{sec:lowStretch}
In this section, we present a simple preconditioning algorithm,
  \texttt{UltraSimple},
  that works by simply adding edges to low-stretch spanning trees.
This algorithm is sufficient to obtain all our results for planar graphs.
For arbitrary graphs, this algorithm might add too many additional edges.
We will show in Section~\ref{sec:ultra} how these extra edges can be removed
  via sparsification.

\subsection{Low-Stretch Trees}\label{ssec:low-stretch}
Low-stretch spanning trees were introduced by Alon, Karp, Peleg and West~\cite{AKPW}.
At present, the construction of spanning trees with the lowest stretch
  is due to Abraham and Neiman~\cite{AbrahamNeiman}, who prove

\begin{theorem}[Low Stretch Spanning Trees]\label{thm:lowStretch}
There exists an $O (m \log n \log \log n)$-time algorithm, \texttt{LowStretch},
  that on input a weighted connected graph $G = (V,E)$,  outputs a spanning
  tree $T$ of $G$ such that
\[
  \stretch{T}{E} 
 \leq c_{AN} \  m \log n \log \log n,
\]
where $m = \sizeof{E}$, for some constant $c_{AN}$.
\end{theorem}

\subsection{Augmenting Low-Stretch Spanning Trees}\label{}

Our procedure for deciding which edges to add to a tree begins by
  decomposing the tree into sub-trees.
In the decomposition, we allow subtrees to overlap at a single vertex,
  or even consist of just a single vertex.
Then, for every pair of subtrees connected by edges of $E$, we add one
  such edge of $E$ to the tree.
The subtrees are specified by the subset of the vertices that they span.

\begin{definition}
Given a tree  $T$ that spans a set of vertices $V$,
a \textit{$T$-decomposition} is a 
  decomposition of $V$ into sets $W_{1}, \dotsc , W_{h}$
  such that $V = \union W_{i}$,
  the graph induced by $T$ on each $W_{i}$ is a tree, possibly
  with just one vertex, 
  and for all $i \not = j$,
  $\sizeof{W_{i} \intersect W_{j}} \leq 1$.

Given an additional set of edges $E$ on $V$,
  a \textit{$(T,E)$-decomposition} is
  a pair
  $(\setof{W_{1}, \dots , W_{h}} , \rho  )$ where
  $\setof{W_{1}, \dots , W_{h}}$ is a $T$-decomposition
  and $\rho  $ is a map 
  that sends each edge of $E$ to a set
  or pair of sets 
  in $\setof{W_{1}, \dots , W_{h}}$
  so that
  for each edge in $(u,v) \in E$,

\begin{dasenumerate}
\item  if
 $\rho (u,v) = \setof{W_{i}}$ then
  $\setof{u,v} \subseteq  W_{i}$, and
\item  
  if $\rho (u,v) = \setof{W_{i}, W_{j}}$, then
  either $u \in W_{i}$ and $v \in W_{j}$, or
  $u \in W_{j}$ and $v \in W_{i}$.
\end{dasenumerate}
\end{definition}

\begin{figure}[h]
\begin{center}
\epsfig{figure=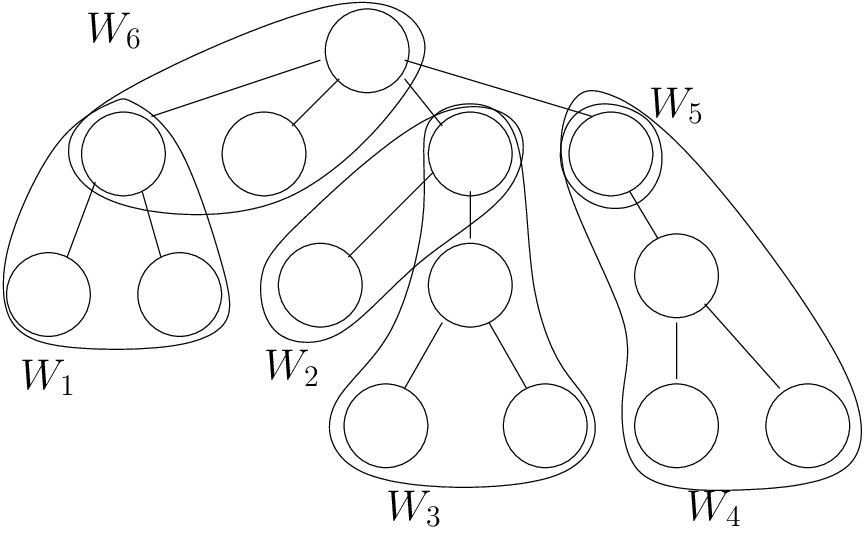,height=1.7in}
\end{center}
\caption{An example of a tree decomposition.  
   Note that sets $W_{1}$ and $W_{6}$ overlap, and that
  set $W_{5}$ is a singleton set and that it overlaps $W_{4}$.
}
\end{figure}

We remark that as the sets $W_{i}$ and $W_{j}$ can overlap,
  it is possible that $\rho (u,v) = \setof{W_{i}, W_{j}}$,
  $u \in W_{i}$ and $v \in W_{i} \intersect W_{j}$.

We use the following tree decomposition theorem
  to show that
  one can always quickly find a
  $T$-decomposition of $E$ with few components 
  in which the sum of stretches of the edges
  attached to each non-singleton component is not too big.
As the theorem holds for any non-negative function $\eta$
  on the edges, not just stretch, we state it in this
  general form.

\begin{theorem}[\texttt{decompose}]\label{thm:decomposeTree}
There exists a linear-time algorithm, which we invoke with the syntax
\[
  (\setof{W_{1}, \ldots , W_{h}}, \rho) = 
   \mathtt{decompose} (T, E, \eta, t),
\]
that on input  a set of edges $E$ on a vertex set $V$, a spanning tree $T$ on $V$,
  a function $\eta : E \rightarrow \Reals{+}$, and
  an integer $1 < t \leq \sum_{e\in E}\eta (e)$,
  outputs a $(T, E)$-decomposition  $(\setof{W_{1}, \ldots , W_{h}}, \rho)$, such that
\begin{dasenumerate}

\item  $h \leq t$,
\item  for all $W_{i}$ such that $\sizeof{W_{i}} > 1$, 
\[
   \sum_{e \in E : W_{i} \in \rho (e) }
   \eta (e)
  \leq 
    \frac{4}{t}
   \sum_{e \in E}
   \eta (e).
\]
\end{dasenumerate}
\end{theorem}
For pseudo-code and a proof of this theorem, see Appendix~\ref{sec:decompose}.
We remark that when $t \geq n$, the algorithm can just construct a singleton
  set for every vertex.

For technical reasons,
  edges with stretch less than 1 can be inconvenient.
So, we define
\begin{equation}\label{eqn:eta}
  \eta (e) = \max (\stretch{T}{e}, 1)
\quad
\text{and}
\quad 
  \eta (E) = \sum_{e \in E} \eta (e).
\end{equation}
The tree $T$ should always be clear from context.

Given a $(T,E)$-decomposition,  $(\setof{W_{1}, \ldots , W_{h}}, \rho)$,
  we define the map 
\[
\sigma : \setof{1,\dotsc ,h} \times \setof{1,\dotsc ,h} \rightarrow
E \union \setof{\mathit{undefined}}
\]
by setting
\begin{equation}\label{eqn:sigma}
\sigma (i,j)
= 
\begin{cases}
\arg \max_{e : \rho (e) = \setof{W_{i}, W_{j}}}
\weight{e} / \eta (e),  & 
\text{if $i \neq j$ and such an $e$ exists}
\\
\mathrm{undefined}  & \text{otherwise.}
\end{cases}
\end{equation}
In the event of a tie, we let $e$ be the lexicographically least edge
 maximizing $\weight{e} / \eta (e)$
  such that $\rho (e) = \setof{W_{i}, W_{j}}$.
Note that $\sigma (i,j)$ is a weighted edge.

The map $\sigma$ tells us which edge from $E$ between $W_{i}$ and $W_{j}$
  to add to $T$.
The following property of $\sigma $, which follows immediately from its definition,
  will be used in our analysis in
  this and the next section.

\begin{proposition}\label{pro:maxbridge}
For every $i,j$ such that $\sigma (i,j)$ is defined and for every $e\in E$ such that 
  $\rho (e) = \setof{W_{i},W_{j}}$,
{\rm
\[
  \frac{\weight{e}}{\eta (e)} \leq  \frac{\weight{\sigma (i,j)}}{\eta (\sigma (i,j))}.
\]
}
\end{proposition}

We can now state the procedure by which we augment
  a spanning tree.

\begin{algbox}
\noindent $F = \texttt{AugmentTree}(T, E, t)$,
\begin{tightlist}
\item []$E$ is set of weighted edges,
\item []$T$ is a spanning tree of the vertices underlying $E$,
\item []$t$ is an integer.

\end{tightlist}
\begin{enumerate}
\item [1.] Compute $\stretch{T}{e}$ for each edge $e \in E$. 

\item [2.]  $(\left(W_{1}, \dotsc , W_{h} \right) , \rho ) = \texttt{decompose} (T,E,\eta ,t)$, where $\eta (e)$ is as defined in \eqref{eqn:eta}.
\item [3.] 
  Set $F$ to be the union of the weighted edges
  $\sigma (i,j)$ over all pairs $1 \leq i < j \leq h$ for which
  $\sigma (i,j)$ is defined,
 where $\sigma (i,j)$ is as defined in \eqref{eqn:sigma}.

\end{enumerate}
\end{algbox}

\begin{algbox}
\noindent $A = \texttt{UltraSimple}(E, t)$
\begin{enumerate}
\item [1.] Set $T = \mathtt{LowStretch} (E)$.
\item [2.] Set $F = \mathtt{AugmentTree} (T, E, t)$.
\item [3.] Set $A = T \union F$.
\end{enumerate}
\end{algbox}

We remark that when $t \geq n$, \texttt{UltraSimple} can just return
  $A = E$.

\begin{theorem}[\texttt{AugmentTree}]\label{thm:augmentTree}
On input a set of weighted edges $E$, a spanning subtree $T$,
  and an integer $1< t \leq \eta (E)$,
  the algorithm \texttt{AugmentTree} runs in time $O (m \log n)$,
  where $m = \sizeof{E}$.
The set of edges $F$ output by the algorithm satisfies
\begin{itemize}
\item [(a)] $F \subseteq E$,
\item [(b)] $\sizeof{F} \leq t^{2}/2$, 
\item [(c)] If $T \subseteq E$, as happens when \texttt{AugmentTree} 
  is called by \texttt{UltraSimple}, then $( T \cup F ) \pleq E$.
\item [(d)] 
\begin{equation}\label{eqn:augmentTree}
 E \pleq 
  \frac{12 \eta (E)}{t}\cdot (T \cup F).
\end{equation}

\end{itemize}
Moreover, if $E$ is planar then $A$ is planar and $\sizeof{F} \leq 3t - 6$.
\end{theorem}

\begin{proof}
In Appendix~\ref{sec:stretch}, we present an algorithm for computing the stretch of
  each edge of $E$ in time $O (m \log n)$.
The remainder of the analysis of the running time is trivial.
Part $(a)$ follows immediately from the statement of the algorithm.
When $T \subseteq E$, $T \cup F \subseteq E$, so part $(c)$ follows as well.

To verify $(b)$, note that the algorithm adds at most one edge to $F$
  for each pair of sets in $W_{1}, \dots , W_{h}$, and there are
  at most ${t \choose 2} \leq t^{2}/2$ such pairs.
If $E$ is planar, then $F$ must be planar as 
  $F$ is a subgraph of $E$.
Moreover, we can use Lemma \ref{lem:planar}
  to show that the graph induced by $E$ on the
  sets $W_{1}, \dots , W_{h}$ is also planar.
Thus,  the number of pairs of these sets connected by edges
  of $E$ is at most the maximum number of edges in a planar 
  graph with $t$ vertices, $3t - 6$.

We now turn to the proof of part $(d)$.
Set 
\begin{equation}\label{eqn:beta}
\beta = 4 \eta (E) /t.
\end{equation}
By Theorem~\ref{thm:decomposeTree}, 
  $\rho$ and $W_{1}, \dotsc , W_{h}$ satisfy
\begin{equation}\label{eqn:sumInW}
  \sum_{e : W_{i} \in \rho (e)}
    \eta (e)
  \leq 
    \beta ,
\quad 
\text{for all $W_{i}$ such that $\sizeof{W_{i}} > 1$.}
\end{equation}

Let $E^{int}_{i}$ denote the set of edges $e$
  with $\rho (e) = \setof{W_{i}}$, and let
  $E^{ext}_{i}$ denote the set of edges $e$
  with $\sizeof{\rho (e)} = 2$
  and $W_{i} \in \rho(e)$.
Let $E^{int} = \union_{i} E^{int}_{i}$
  and $E^{ext} = \union_{i} E^{ext}_{i}$.
Also, let $T_{i}$ denote the tree formed by the edges
  of $T$ inside the set $W_{i}$.
Note that 
  when $\sizeof{W_{i}} = 1$, $T_{i}$ and $E^{int}_{i}$
  are empty.

We will begin by proving that when $\sizeof{W_{i}} > 1$,
\begin{equation}\label{eqn:augmentTreeInt}
  E^{int}_{i} \pleq \left(\sum_{e \in E^{int}_{i}} \eta (e) \right) T_{i},
\end{equation}
from which it follows that
\begin{equation}\label{eqn:augmentTreeInt2}
  E^{int}
 \pleq 
 \sum_{i : \sizeof{W_{i}} > 1}  \left(\sum_{e \in E^{int}_{i}} \eta (e) \right) T_{i}.
\end{equation}

For any edge  $e \in E^{int}_{i}$, 
  the path in $T$ between the endpoints of $e$ lies entirely in $T_{i}$.
So, by \eqref{eqn:edgeStretch} we have
\[
  e \pleq \stretch{T}{e}\cdot  T_{i} \pleq \eta (e) \cdot  T_{i}.
\]
Inequality~\eqref{eqn:augmentTreeInt} now follows by summing
  over the edges $e \in E^{int}_{i}$.

We now define the map $\tau : E \rightarrow E \union \setof{\mathit{undefined}}$ by
\begin{equation}\label{eqn:tau}
\tau (e) = 
\begin{cases}
\sigma (i,j),  &  
\text{if $\sizeof{\rho (e)} = 2$, where $\rho (e) = \setof{W_{i}, W_{j}}$, and}        \\
\mathit{undefined}  & \text{otherwise}.
\end{cases}
\end{equation}

To handle the edges bridging components, 
  we prove that for each edge $e$ with $\rho (e) = (W_{i}, W_{j})$,
\begin{equation}\label{eqn:augmentTreeExt}
 e \pleq 3 \eta (e) (T_{i} + T_{j}) + 
  3 \frac{\weight{e}}{ \weight{\tau (e)}} \cdot \tau (e) 
\end{equation}
Let $e = w \edg{u,v}$ be such an edge, with $u \in W_{i}$ and $v \in W_{j}$.
Let $\tau (e) = z \edg{x,y}$, with $x \in W_{i}$ and $y \in W_{j}$.
Let $t_{i}$ denote the last vertex in $T_{i}$ on the path in $T$
  from $u$ to $v$ (see Figure~\ref{fig:treeLabeled}).
If $T_{i}$ is empty, $t_{i} = u$.
Note that $t_{i}$ is also the last vertex in $T_{i}$ on the path in $T$
  from $x$ to $y$.
Define $t_{j}$ similarly.
\begin{figure}[h]
\begin{center}
\epsfig{figure=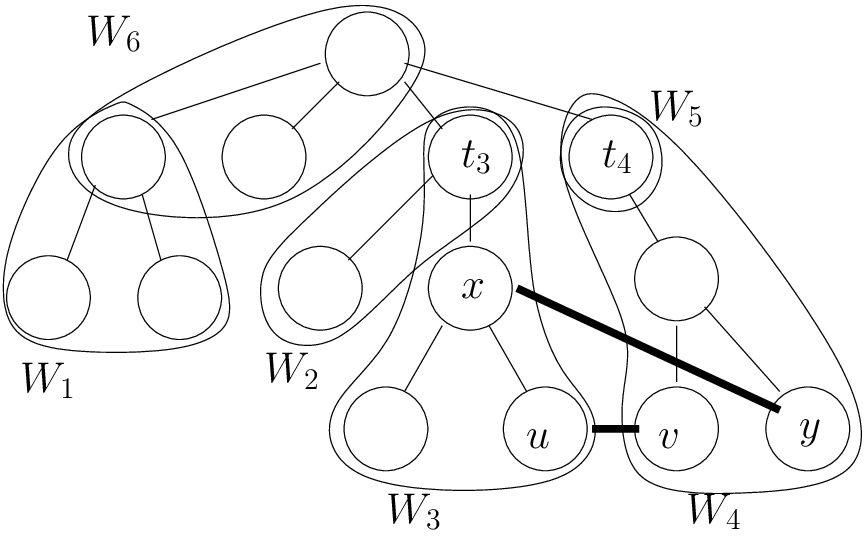,height=1.7in}
\end{center}
\caption{In this example, $e = w \edg{u,v}$ and $\tau (e) = z \edg{x,y}$.}
\label{fig:treeLabeled}
\end{figure}
As $T_{i}(u,x) \subseteq T_{i}(u,t_{i}) \union T_{i}(t_{i},x)$,
 the tree $T_{i}$ contains a path from $u$ to $x$ of resistance at most
\[
\res{T_{i}(u,t_{i})}
+
\res{T_{i}(t_{i},x)},
\]
and the tree $T_{j}$ contains a path from $y$ to $v$ of resistance at most
\[
\res{T_{j}(y,t_{j})}
+
\res{T_{j}(t_{j},v)}.
\]
Furthermore, as $T_{i}(u, t_{i}) + T_{j}(t_{j}, v) \subseteq T(u,v)$
  and $T_{i}(t_{i}, x) + T_{j}(y, t_{j}) \subseteq T(x,y)$,
  the sum of the resistances of the paths from $u$ to $x$ in $T_{i}$
  and from $y$ to $v$ in $T_{j}$ is at most
\begin{align*}
\res{T(u,v)} + \res{T(x,y)}
& =
  \stretch{T}{e}/w + \stretch{T}{\tau (e)}/z
\\
& \leq 
  \eta (e)/w + \eta (\tau (e))/z
\\
& \leq 
 2 \eta (e) / w,
\end{align*}
where the last inequality follows from Proposition \ref{pro:maxbridge}.
Thus, the graph
\[
  3 \eta (e) (T_{i} + T_{j})
  + 3 w \edg{x,y}
=
  3 \eta (e) (T_{i} + T_{j})
+  3 \frac{\weight{e}}{ \weight{\tau (e)}} \cdot \tau (e) 
\]
contains a path from $u$ to $v$  
  of resistance at most
\[
  \frac{2}{3} \frac{1}{w} + \frac{1}{3} \frac{1}{w}
  = \frac{1}{w},
\]
which by Lemma~\ref{lem:pathResistance} implies
  (\ref{eqn:augmentTreeExt}).

We will now sum \eqref{eqn:augmentTreeExt} over 
  every edge $e \in E_{i}^{ext}$ for every $i$, observing that
  this counts every edge in $E^{ext}$ twice.
\begin{align}
E^{ext}
& =
(1/2) \sum_{i} \sum_{e \in E^{ext}_{i}} e
\notag \\
& \pleq
\sum_{i} 
  \sum_{e \in E^{ext}_{i}} 
3  \eta (e) T_{i} 
 + 
(1/2)
\sum_{i} 
  \sum_{e \in E^{ext}_{i}} 
3 \frac{\weight{e}}{ \weight{\tau (e)}} \cdot \tau (e) 
\notag \\
& =
3 \sum_{i} 
  \left( 
   \sum_{e \in E^{ext}_{i}}  \eta (e)
  \right)
T_{i}  
 + 
3  \sum_{e \in E^{ext}} 
 \frac{\weight{e}}{ \weight{\tau (e)}} \cdot \tau (e)
\notag \\
& =
3 \sum_{i: \sizeof{W_{i}} > 1} 
  \left( 
   \sum_{e \in E^{ext}_{i}}  \eta (e)
  \right)
T_{i}  
 + 
3  \sum_{e \in E^{ext}} 
 \frac{\weight{e}}{ \weight{\tau (e)}} \cdot \tau (e),
\label{eqn:augmentTree1}
\end{align}
as $T_{i}$ is empty when $\sizeof{W_{i}} = 1$.

We will now upper bound the right-hand side of \eqref{eqn:augmentTree1}.
To handle boundary cases, we divide $E^{ext}$ into two sets.
We let $E^{ext}_{single}$ consist of those
  $e \in E^{ext}$ for which both sets in $\rho (e)$ have size $1$.
We let $E^{ext}_{general} = E^{ext} - E^{ext}_{single}$ contain
  the rest of the edges in $E^{ext}$.
For $e \in E^{ext}_{single}$, $\tau (e) = e$, while for
  $e \in E^{ext}_{general}$, $\tau (e)\in E^{ext}_{general}$.

For $E^{ext}_{single}$, we have
\[
  \sum_{e \in E^{ext}_{single}} 
   \frac{\weight{e}}{ \weight{\tau (e)}} \cdot \tau (e)
=
  \sum_{e \in E^{ext}_{single}}  \tau (e)
=
   E^{ext}_{single}.
\]

To evaluate the sum over the edges $e \in E^{ext}_{general}$,
  consider any $f \in E^{ext}_{general}$ in the image of $\tau $.
Let $i$ be such that   $f \in E^{ext}_{i}$ and $\sizeof{W_{i}} > 1$.
Then, for every  $e$ such that
  $\tau (e) = f$, we have $e \in E^{ext}_{i}$.
So, by Proposition \ref{pro:maxbridge},
\begin{multline}\label{eqn:augmentTreeSlack}
\sum_{
\substack{e \in  E^{ext}\\
\tau (e) = f}
} 
\frac{\weight{e} }{ \weight{\tau (e)}}
=
\sum_{
\substack{e \in  E^{ext}_{i}\\
\tau (e) = f}
} 
\frac{\weight{e} }{ \weight{\tau (e)}}
\\
\leq 
\sum_{
  e \in  E^{ext}_{i}} 
\frac{\weight{e} }{ \weight{\tau (e)}}
\leq 
\sum_{e \in  E^{ext}_{i}} 
\frac{ \eta (e) }{ \eta (\tau (e))}
\leq 
\sum_{e \in  E^{ext}_{i}} 
 \eta (e) 
\leq 
  \beta .
\end{multline}
Thus,
\[
  \sum_{e \in E^{ext}} 
 \frac{\weight{e}}{ \weight{\tau (e)}} \cdot \tau (e)
\pleq
 E^{ext}_{single} + 
  \sum_{
\substack{f \in \mathrm{image} (\tau )\\
f \in E^{ext}_{general}}} 
 \beta  \cdot f
\pleq 
 \beta  \cdot F.
\]
Plugging this last inequality into \eqref{eqn:augmentTree1},
we obtain
\[
E^{ext} \pleq 
3 \sum_{i : \sizeof{W_{i}} > 1} 
  \left( 
   \sum_{e \in E^{ext}_{i}}  \eta (e)
  \right)
T_{i}  
 + 
3 \beta  \cdot F.
\]
Applying
  \eqref{eqn:augmentTreeInt2} and then \eqref{eqn:sumInW}, we compute
\[
E 
=
E^{ext} + E^{int}
\pleq
3 \sum_{i : \sizeof{W_{i}} > 1} T_{i}  
  \left( 
   \sum_{e \in E^{int}_{i}}  \eta (e)
 +
   \sum_{e \in E^{ext}_{i}}  \eta (e)
  \right)
 + 
3 \beta  \cdot F
\pleq 3 \beta  \cdot  (T \cup F),
\]
which by  \eqref{eqn:beta} implies the lemma.
\end{proof}

We now discuss a source of slack in Theorem~\ref{thm:augmentTree}.
This is the motivation for the construction of ultra-sparsifiers
  in the next section.

  In the proof of Theorem~\ref{thm:augmentTree}, we assume in the worst
  case that the tree decomposition could result in 
  each tree $T_{i}$ being connected to $t-1$ other trees, for a total
  of $t (t-1)/2$ extra edges.
  Most of these edges seem barely necessary, as they could be included at a
  small fraction of their original weight. 
  To see why, consider the crude estimate at the end of inequality   
  \eqref{eqn:augmentTreeSlack}.  We upper bound the multiplier of one bridge
  edge $f$ from $T_{i}$,
\[
\sum_{
\substack{e \in  E^{ext}_{i}\\
\tau (e) = f}
} 
\frac{\weight{e} }{ \weight{\tau (e)}},
\]
by the sum of the multipliers of all bridge edges from $T_{i}$,
\[
\sum_{
  e \in  E^{ext}_{i}} 
\frac{\weight{e} }{ \weight{\tau (e)}}.
\]
The extent to which this upper bound is loose is the factor by which we could
  decrease the weight of the edge $f$ in the preconditioner.

  While we do not know how to accelerate our algorithms by decreasing the weights with
  which we include edges, we are able to use sparsifiers to trade many low-weight
  edges for a few edges of higher weight.
  This is how we reduce the number of edges we add to the spanning tree to
  $O (t \log^{c_{2} + 5} n)$.

\section{Ultra-Sparsifiers}\label{sec:ultra}
We begin our construction of ultra-sparsifiers by building
  ultra-sparsifiers for the special case in which our graph
  has a distinguished vertex $r$ and a low-stretch spanning tree $T$
  with the property that for every edge $e \in E - T$, the path in
  $T$ connecting the endpoints of $e$ goes through $r$.
In this case, we will call $r$ the \textit{root} of the tree.
All of the complications of ultra-sparsification will be handled
  in this construction.
The general construction will follow simply
  by using tree splitters to choose the roots and decompose
  the input graph.

\begin{figure}[htop]
\begin{algbox}
\noindent $E_{s} = \texttt{RootedUltraSparsify}(E, T, r, t, p)$\\

\noindent Condition: 
  for all $e \in E$, $r \in T (e)$. The parameter $t$ is a positive integer at
  most $\ceiling{\eta (E)}$.
\begin{enumerate}

\item Compute $\stretch{T}{e}$ and $\eta (e)$ for each edge $e \in E$,
  where $\eta$ is as defined in \eqref{eqn:eta}.

\item If $t \geq \sizeof{E}$, return $E_{s} = E$.

\item Set $(\setof{W_{1}, \dotsc , W_{h}} , \rho ) =
  \texttt{decompose} (T,E,\eta ,t)$.

\item Compute $\sigma$, as given by \eqref{eqn:sigma}, everywhere it is defined.

\item For every $(i,j)$ such that $\sigma (i,j)$ is defined,
  set
\begin{equation}\label{eqn:omega}
\omega (i,j) =   
  \sum_{e \in E : \rho (e) = \setof{W_{i},W_{j}}}
    \weight{e}
\quad 
\text{and}
\quad
\psi (i,j) = \omega (i,j) / \weight{\sigma (i,j)}.
\end{equation}

\item 
Set $F = \setof{\psi  (i,j) \sigma (i,j) : \sigma (i,j) \text{ is defined}}.$
\label{step:RUSdefF}

\item For each $f = \psi  (i,j) \sigma (i,j) \in F$, set
\begin{equation}\label{eqn:phi}
\phi (f) = \max (\psi (i,j), \stretch{T}{f}).
\end{equation}

\item For $b \in \setof{1,\dotsc ,\ceiling{\log_{2} \eta (E)}}$:

\begin{enumerate}
\item 
Set
$
F^{b}
= \begin{cases}
\setof{f \in F : \phi(f) \in [1,2]}  & \text{if $b = 1$}
\\
\setof{f \in F : \phi(f) \in (2^{b-1}, 2^{b}]}  & \text{otherwise}
\end{cases}
$
\item Let $H^{b}$ be the set of edges on vertex set $\setof{1,\dotsc ,h}$
  defined by
\[
H^{b} = \setof{\omega (i,j) \edg{i,j} : \psi (i,j) \sigma (i,j) \in F^{b}}.
\]

\item Set $H^{b}_{s} = \mathtt{Sparsify2} (H^{b}, p)$.

\item Set
\[
E^{b}_{s} = \setof{\sigma (i,j) : \exists w \text{ such that } w\edg{i,j}\in H^{b}_{s}}.
\]
\end{enumerate}
\item Set $E_{s} = \union_{b} E^{b}_{s}$.
\end{enumerate}
\end{algbox}
\end{figure}

The algorithm \texttt{RootedUltraSparsify} begins by computing
  the same set of edges $\sigma (i,j)$, as was computed by \texttt{UltraSimple}.
However, when \texttt{RootedUltraSparsify} puts one of these edges into the set $F$,
  it gives it a different weight: $\omega (i,j)$.
For technical reasons, the set $F$ is decomposed into subsets $F^{b}$
  according to the quantities $\phi (f)$, which will play 
  a role in the analysis of \texttt{RootedUltraSparsify}
  analogous to the  role played by 
  $\eta (e)$ in the analysis of \texttt{UltraSimple}.
Each set of edges $F^{b}$ is sparsified, and the union of 
  the edges of $E$ that appear in the resulting sparsifiers
  are returned by the algorithm.
The edges in $F^{b}$ cannot necessarily be sparsified directly,
  as they might all have different endpoints.
Instead, $F^{b}$ is first projected to a graph $H^{b}$ on vertex set
  $\setof{1,\dotsc , h}$.
After a sparsifier $H^{b}_{s}$ of $H^{b}$ is computed, it is lifted back
  to the original graph to form $E^{b}_{s}$.
Note that the graph $E_{s}$ returned by \texttt{RootedUltraSparsify} is
  a subgraph of $E$, with the same edge weights.

We now prove that $F = \union_{b=1}^{\ceiling{\log_{2} \eta (E)}} F^{b}$.
Our proof will use the function $\eta$, which we recall was defined
  in~\eqref{eqn:eta} and which was used to define the map $\sigma$.

\begin{lemma}\label{lem:phi}
For $\phi$ as defined in \eqref{eqn:phi},
  for every $f = \psi (i,j) \sigma (i,j)\in F$,
\begin{equation}\label{eqn:maxphi}
1 \leq  \psi (i,j) \leq \phi (f) \leq \eta (E).
\end{equation}
\end{lemma}
\begin{proof}
Recall from the definitions of $\phi$ and $\psi$ that
\[
\phi (f) \geq \psi (i,j)
= 
\frac{  \sum_{e \in E : \rho (e) = \setof{W_{i},W_{j}}}
    \weight{e}}{\weight{\sigma (i,j)}}.
\]
By definition $\sigma (i,j)$ is an edge in $E$ satisfying
  $\rho (\sigma (i,j)) = \setof{W_{i},W_{j}}$;
so, the right-hand side of the last expression is at least 1.

To prove the upper bound on $\phi (f)$, first apply Proposition~\ref{pro:maxbridge} 
  to show that 
\begin{equation*}
\psi (i,j) = 
\frac{\sum_{e \in E : \rho (e) = \setof{W_{i},W_{j}}} \weight{e}}
     {\weight {\sigma (i,j)} }
\leq
\frac{\sum_{e \in E : \rho (e) = \setof{W_{i},W_{j}}} \eta (e)}
     {\eta (\sigma (i,j))}
 \leq \eta (E),
\end{equation*}
as $\eta$ is always at least $1$.
Similarly,
\begin{multline*}
\stretch{T}{f} = 
\frac{\omega (i,j)}{\weight{\sigma (i,j)}} \stretch{T}{\sigma (i,j)}
=
\frac{\stretch{T}{\sigma (i,j)}}
     {\weight {\sigma (i,j)}}
\left(\sum_{e \in E : \rho (e) = \setof{W_{i},W_{j}}}
    \weight{e} \right)
\\
\leq 
\frac{\eta (\sigma (i,j))}
     {\weight {\sigma (i,j)}}
\left(\sum_{e \in E : \rho (e) = \setof{W_{i},W_{j}}}
    \weight{e} \right)
\leq \sum_{e \in E : \rho (e) = \setof{W_{i},W_{j}}}
   \eta (e) \leq \eta (E),
\end{multline*}
where the second-to-last inequality follows from Proposition~\ref{pro:maxbridge}.
\end{proof}

It will be convenient for us to extend the domain of $\rho$ to $F$ 
  by setting $\rho (f) = \rho (e)$ where $e \in E$ has the same vertices as $f$.
That is, when there exists
  $\gamma \in \Reals{+}$ such that $f = \gamma e$.
Define 
\[
\beta = 4 \eta (E) / t.
\]
Our analysis of \texttt{RootedUltraSparsify} will exploit the 
  inequalities contained in the following two lemmas.
\begin{lemma}\label{lem:sumF}
For every $i$ for which $\sizeof{W_{i}} > 1$,
\[
  \sum_{f \in F : W_{i} \in \rho (f)} \stretch{T}{f} \leq 
  \beta .
\]
\end{lemma}
\begin{proof}
Consider any $f \in F$, and let $f = \psi (i,j) \sigma (i,j)$.
Note that the weight of $f$ is $\omega (i,j)$, and recall
  that $\stretch{T}{f} \leq \eta (f)$.
We first show that
\[
\sum_{e : \tau (e) = \sigma (i,j)} \eta (e)
\geq \eta (f)  .
\]
By Proposition \ref{pro:maxbridge}, and the definition of $\tau$ in \eqref{eqn:tau} 
\begin{align*}
\sum_{e : \tau (e) = \sigma (i,j)} \eta (e)
& \geq 
\frac{\eta (\sigma (i,j))}{\weight{\sigma (i,j)}}
\sum_{e : \tau (e) = \sigma (i,j)} \weight{e}
\\
& =
\frac{\eta (\sigma (i,j))}{\weight{\sigma (i,j)}}
\weight{f} 
\\
& =
\max  \left( 
  \frac{\weight{f}}{\weight{\sigma (i,j)}} ,
  \frac{\stretch{T}{\sigma (i,j)}}{\weight{\sigma (i,j)}} \weight{f}
  \right)
\\
& =
\max  \left( 
  \psi (i,j),
  \stretch{T}{f}
  \right)
\\
& =
\max  \left( 
  \phi (f),
  \stretch{T}{f}
  \right) \quad \text{(by \eqref{eqn:phi})}
\\
& \geq 
\max  \left( 
  1,
  \stretch{T}{f}
  \right) \quad \text{(by \eqref{eqn:maxphi})}
\\
& =
\eta (f).
\end{align*}
We then have
\[
  \sum_{e \in E : W_{i} \in \rho (e)} \eta (e)
\geq 
  \sum_{f \in F : W_{i} \in \rho (f)} \eta (f).
\]
The lemma now follows from the upper bound 
  of $4 \eta (E) / t$
  imposed on the
  left-hand term by Theorem~\ref{thm:decomposeTree}.
\end{proof}

\begin{lemma}\label{lem:phi2}
For every $i$ for which $\sizeof{W_{i}} > 1$,
\begin{equation}\label{eqn:sumphi}
  \sum_{f \in F : W_{i} \in \rho (f)}
  \phi (f)
\leq 
 2 \beta .
\end{equation}
\end{lemma}
\begin{proof}
For an edge $f \in F$, let $\psi (f)$ equal
  $\psi (i,j)$ where $f = \psi (i,j) \sigma (i,j)$.
With this notation, we may compute
\begin{align*}
  \sum_{f \in F : W_{i} \in \rho (f)}
  \phi (f)
& \leq 
  \sum_{f \in F : W_{i} \in \rho (f)}
   \stretch{T}{f}
 +
  \sum_{f \in F : W_{i} \in \rho (f)}
    \psi (f)\\
& \leq 
  \sum_{f \in F : W_{i} \in \rho (f)}
   \eta (f)
 +
  \sum_{f \in F : W_{i} \in \rho (f)}
    \psi (f)\\
& \leq 
\beta 
 +
  \sum_{f \in F : W_{i} \in \rho (f)}
    \psi (f),
\end{align*}
by Lemma~\ref{lem:sumF}.
We now bound the right-hand term as
  in the proof of inequality~\eqref{eqn:augmentTreeSlack}:
\[
\sum_{f \in F: W_{i} \in \rho (f)}
  \psi (f)
=
\sum_{e \in E^{ext}_{i}}
\frac{\weight{e}}{\weight{\tau (e)}}
\\
\leq 
\sum_{e \in E^{ext}_{i}}
\frac{\eta (e)}{\eta (\tau (e))}
\leq 
\sum_{e \in E^{ext}_{i}}
\eta (e) 
\leq
\beta ,
\]
by our choice of $\beta $ and Theorem~\ref{thm:decomposeTree}.
\end{proof}

\begin{lemma}[\texttt{RootedUltraSparsify}]\label{lem:rootedUltra}
Let $T$ be a spanning tree on a vertex set $V$, and let
  $E$ be a nonempty set of edges on $V$ for which there exists an 
  $r \in V$ be such that
  for all $e \in E$, $r \in T (e)$.
For $p >0$ and $t$ a positive integer at most $\ceiling{\eta (E)}$,
  let $E_{s}$ be the
  graph returned by $\texttt{RootedUltraSparsify}(E, T, r, t, p)$.
The graph $E_{s}$ is a subgraph of $E$, and
  with probability at least $1 - \ceiling{\log_{2} \eta (E)} p$, 
\begin{equation}\label{eqn:rootedUltraSize}
\sizeof{E_{s}} \leq 
c_{1} \log^{c_{2}} (n/p) \max (1,\ceiling{\log_{2} \eta (E)}) t,
\end{equation}
and
\begin{equation}\label{eqn:rootedUltraLeq}
 E \pleq 
  \left(3 \beta  + 126 \beta  \max (1,\log_{2} \eta (E)) \right) \cdot T + 120 \beta  \cdot E_{s},
\end{equation}
where $\beta  = 4  \eta (E) /t$.
\end{lemma}
\begin{proof}
We first dispense with the case in which the algorithm terminates at line 2.
If $t \geq m$, then both \eqref{eqn:rootedUltraSize} and \eqref{eqn:rootedUltraLeq}
  are trivially satisfied by setting $E_{s} = E$, as
  $\beta \geq 2$.

By Theorem~\ref{thm:sparsifier} each graph $H^{b}_{s}$
  computed by \texttt{Sparsify2} is a 
  $c_{1} \log^{c_{2}} (n/p)$-sparsifier of $H^{b}$ according to
  Definition~\ref{def:sparsifier}
  with probability at least $1-p$.
As there are at most $\ceiling{\log_{2}\eta (E)}$ such graphs $H^{b}$,
  this happens for all of these graphs with probability at least
  $1 - \ceiling{\log_{2} \eta (E)} p$.
For the remainder of the proof, we will assume that each graph
  $H^{b}_{s}$
  is a 
  $c_{1} \log^{c_{2}} (n/p)$-sparsifier of $H^{b}$.
Recalling that $h \leq t$, the bound on the number of edges in $E_{s}$
  is immediate.

Our proof of \eqref{eqn:rootedUltraLeq} will go through an analysis
  of intermediate graphs.
As some of these could be multi-graphs, we will find it convenient
  to write them as sums of edges.

To define these intermediate graphs, 
  let $r_{i}$ be the vertex in $W_{i}$ that is closest
  to $r$ in $T$.
As in Section~\ref{sec:lowStretch}, 
 let $T_{i}$ denote the edges of the 
 subtree of $T$ with vertex set $W_{i}$.
We will view $r_{i}$ as the root of tree $T_{i}$.
Note that if $\sizeof{W_{i}} = 1$,
  then $W_{i} = \setof{r_{i}}$ and $T_{i}$ is empty.
As distinct sets $W_{i}$ and $W_{j}$ can overlap in at most one
  vertex, $\sum_{i} T_{i} \leq T$.
We will exploit the fact that for each $e \in E$
  with $\rho (e) = \setof{W_{i}, W_{j}}$,
  the path $T (e)$ contains both $r_{i}$ and $r_{j}$,
which follows from the condition $r \in T (e)$.

We now define the edge set $D^{b}$, which is a projection
  of $H^{b}$ to the vertex set $r_{1}, \dotsc , r_{h}$,
  and $D^{b}_{s}$, which is an analogous projection of the sparsifier 
  $H^{b}_{s}$.
We set
\[
D^{b} = \sum_{(i,j) : \psi (i,j)\sigma (i,j) \in F^{b}} 
      \omega (i,j) \edg{r_{i}, r_{j}}
\]
and
\[
D^{b}_{s} = \sum_{w \edg{i,j} \in H^{b}_{s}} w \edg{r_{i}, r_{j}}.
\]
As the sets $W_{i}$ and $W_{j}$ are allowed to overlap slightly,
  it could be the case that some $r_{i} = r_{j}$ for $i \not = j$.
In this case, $D^{b}$ would not be isomorphic to $H^{b}$.

Set
\[
  F^{b}_{s} = \setof{\gamma  \psi  (i,j) \sigma (i,j) : 
\text{$\exists \gamma$ and $(i,j)$ so that $\gamma \omega (i,j) \edg{i,j} \in H^{b}_{s}$}}.
\]
The edge set $H^{b}$ can be viewed as a projection of the edge set $F^{b}$ to the vertex
  set $\setof{1,\dotsc , h}$, and the edge set $F^{b}_{s}$ can be viewed as a lift
  of $H^{b}_{s}$ back into a reweighted subgraph of $F^{b}$.

We will prove the following inequalities
\begin{align}
E & \pleq 3 \beta  \cdot  T + 3 F \label{eqn:E0toTF}\\
F^{b} & \pleq 2 \beta  \cdot T + 2 D^{b} \label{eqn:FtoD}\\
D^{b} & \pleq (5/4) D^{b}_{s} \label{eqn:DtoDs}\\
D^{b}_{s} & \pleq 16 \beta  \cdot T + 2 F^{b}_{s} \label{eqn:DstoTFs}\\
F^{b}_{s} & \pleq 8 \beta  \cdot E^{b}_{s} \label{eqn:FstoEs}
\end{align}
Inequality~\eqref{eqn:rootedUltraLeq} in the statement of the lemma follows 
  from these inequalities and $F = \sum_{b} F^{b}$.

To prove inequality~\eqref{eqn:E0toTF}, we exploit
  the proof of Theorem~\ref{thm:augmentTree}.
The edges $F$ constructed in \texttt{RootedUltraSparsify}
  are the same as those chosen by \texttt{UltraSimple},
  except that they are reweighted by the function $\psi$.
If we follow the proof of inequality~\eqref{eqn:augmentTree} in
  Theorem~\ref{thm:augmentTree}, but
  neglect to apply inequality \eqref{eqn:augmentTreeSlack},
  we obtain
\[
E \pleq 3 \beta  \cdot T + 3  \sum_{e \in E^{ext}} 
 \frac{\weight{e}}{ \weight{\tau (e)}} \cdot \tau (e)
= 
3 \beta  \cdot T + 3 F.
\]

To prove inequality~\eqref{eqn:FtoD},
  consider any edge $w \edg{u,v} = f \in F^{b}$.
Assume $\rho (f) = \setof{W_{i}, W_{j}}$,  $u \in W_{i}$ and $v \in W_{j}$.
We will now show that
\begin{equation}\label{eqn:FtoDstep}
f \pleq 2 \stretch{T}{f} (T_{i} +  T_{j}) 
  + 2 w \edg{r_{i}, r_{j}}.
\end{equation}
As the path from $u$ to $v$ in $T$ contains both $r_{i}$ and $r_{j}$,
\[
\res{T (u, r_{i})}
+
\res{T (r_{j}, v)}
\leq 
\res{T (u,v)}
=
\stretch{T}{f} / w.
\]
Thus, the resistance of the path
\[
2 \stretch{T}{f} T (u, r_{i})
+
2 w \edg{r_{i}, r_{j}}
+
2 \stretch{T}{f} T (r_{j}, v)
\]
is at most $1/w$, and so
 Lemma~\ref{lem:pathResistance} implies that
\[
f
\pleq 
2 \stretch{T}{f} T (u, r_{i})
+
2 w \edg{r_{i}, r_{j}}
+
2 \stretch{T}{f} T (r_{j}, v),
\]
which in turn implies \eqref{eqn:FtoDstep}.
Summing \eqref{eqn:FtoDstep} over all $f \in F^{b}$
  yields
\begin{align*}
  F^{b}
&  \pleq
 2  \sum_{i}
  \left(
   \sum_{f \in F : W_{i} \in \rho (f)}
   \stretch{T}{f}  
   \right)
  T_{i}
+
  2 D^{b}
\\
  F^{b}
&  \pleq
 2  \sum_{i : \sizeof{W_{i}} > 1}
  \left(
   \sum_{f \in F : W_{i} \in \rho (f)}
   \stretch{T}{f}  
   \right)
  T_{i}
+
  2 D^{b} &
\text{as $T_{i}$ is empty when $\sizeof{W_{i}} = 1$}
\\ 
& \pleq 
 2  \sum_{i} \beta  \cdot T_{i}
+
  2 D^{b},
& \text{by Lemma~\ref{lem:sumF}}
\\
&\pleq 
 2 \beta  \cdot T + 2 D^{b}.
\end{align*}

We now prove inequality \eqref{eqn:DstoTFs}, as it uses similar techniques.
Let $f_{s} = w \edg{u,v} \in F^{b}_{s}$.
Then, there exist  $\gamma$ and $(i,j)$ so  that
  $\gamma \omega (i,j) \edg{i,j} \in H^{b}_{s}$,
  $u \in W_{i}$, and $v \in W_{j}$.
Set $\gamma (f_{s})$ to be this multiplier $\gamma$.
By part $(c)$ of Definition~\ref{def:sparsifier},
  we must have
  $\omega (i,j) \edg{i,j} \in H^{b}$ and
  $\psi (i,j) \sigma (i,j) \in F^{b}$.
Let $f = \psi (i,j) \sigma (i,j)$.
Note that $f_{s} = \gamma (f_{s}) f$.
The sum of the resistances of the paths from
  $r_{i}$ to $u$ in $T_{i}$
  and from $v$ to $r_{j}$ in $T_{j}$
 is 
\[
\res{T (r_{i}, u)}
+
\res{T (v, r_{j})}
\leq 
\res{T (u,v)}
=
\stretch{T}{f} / \omega (i,j),
\]
as $\weight{f} = \omega (i,j)$.
Thus, the resistance of the path
\[
2 \stretch{T}{f} T (r_{i}, u)
+
2  f
+
2 \stretch{T}{f} T (v, r_{j})
\]
is at most $1/\omega (i,j)$, and so
 Lemma~\ref{lem:pathResistance} implies that
\[
\omega (i,j) \edg{r_{i}, r_{j}}
\pleq 
2 \stretch{T}{f} (T_{i} + T_{j}) + 2 f,
\]
and
\begin{align*}
\gamma (f_{s}) \omega (i,j) \edg{r_{i}, r_{j}}
& \pleq 
2 \gamma (f_{s}) \stretch{T}{f} (T_{i} + T_{j}) + 2 f_{s}
\\
& \pleq 
2 \gamma (f_{s}) \phi (f) (T_{i} + T_{j}) + 2 f_{s}
& \text{(by \eqref{eqn:phi})}
\\
& \pleq 
2^{b+1} \gamma (f_{s})  (T_{i} + T_{j}) + 2 f_{s}
& \text{(by $f \in F^{b}$).}
\end{align*}
Summing this inequality over all $f_{s} \in F^{b}_{s}$, we obtain
\[
D^{b}_{s} \pleq 
\sum_{i} 
\left(2^{b+1} \sum_{f_{s} \in F^{b}_{s}: W_{i} \in \rho (f_{s})} \gamma (f_{s})
\right) T_{i}
+
2 F^{b}_{s}.
\]
For all $i$ such that $\sizeof{W_{i}} > 1$,
\begin{align}
  \sum_{f_{s} \in F^{b}_{s} : W_{i} \in \rho (f_{s})} \gamma (f_{s})
& \leq 
  2 \sizeof{\setof{f \in F^{b} : W_{i} \in \rho (f)  }}
& \text{(part $(d)$ of Definition~\ref{def:sparsifier})}
\notag
\\
& \leq
  2 \sum_{f \in F^{b} : W_{i} \in \rho (f)}
   \phi (f) / 2^{b-1}
\notag
\\
& \leq 
  4 \beta  / 2^{b-1}
& \text{(by Lemma~\ref{lem:phi2})}
\notag
\\
&
= \beta  / 2^{b-3}.\label{eqn:sumfs}
\end{align}
So,
\[
D^{b}_{s} \pleq 
\sum_{i} 16 \beta  \cdot T_{i}
+
2 F^{b}_{s}
\pleq 
16 \beta  \cdot T
+ 2 F^{b}_{s}.
\]

To prove inequality~\eqref{eqn:FstoEs},
  let $f_{s}$ be any edge in $F_{s}$,
  let $f$ be the edge in $F$ such that
  $f_{s} = \gamma (f_{s}) f$,
  and let $\sigma (i,j)$ be the edge such that
  $f_{s} = \gamma (f_{s}) \psi (i,j) \sigma (i,j)$.
It suffices to show that
\begin{equation}\label{eqn:FstoEsPart}
\weight{f_{s}}  \leq 8 \beta  \  \weight{\sigma (i,j)}.
\end{equation}
Set $b$ so that $f \in F^{b}$.
By \eqref{eqn:sumfs},
\[
 \gamma (f_{s}) \leq \beta  / 2^{b-3} 
\leq 
  8 \beta  / \phi (f)
=
 8 \beta  / \max (\psi (i,j), \stretch{T}{f})
\leq 
  8 \beta  / \psi (i,j).
\]
As
 $\weight{f_{s}} = \gamma (f_{s}) \psi (i,j) \weight{\sigma (i,j)}$,
 inequality~\eqref{eqn:FstoEsPart} follows.

It remains to prove inequality~\eqref{eqn:DtoDs}.
The only reason this inequality is not immediate from part $(a)$
  of Definition~\ref{def:sparsifier} is that 
  we may have $r_{i} = r_{j}$ for some $i \not = j$.
Let $R = \setof{r_{1}, \dotsc , r_{h}}$
  and $S = \setof{1,\dotsc ,h}$,
Define the map $\pi : \Reals{R} \rightarrow \Reals{S}$
  by $\pi (x)_{i} = x_{r_{i}}$.
We then have for all $x \in \Reals{R}$
\[
x^{T} L_{D^{b}} x = \pi (x)^{T} L_{H^{b}} \pi (x)
\quad 
\text{ and }
\quad 
x^{T} L_{D^{b}_{s}} x = \pi (x)^{T} L_{H^{b}_{s}} \pi (x);
\]
so,
\[
  x^{T} L_{D^{b}} x
=
 \pi (x)^{T} L_{H^{b}} \pi (x)
\leq
(5/4) \pi (x)^{T} L_{H^{b}_{s}} \pi (x)
=
(5/4) x^{T} L_{D^{b}_{s}} x.
\]
\end{proof}

The algorithm \texttt{UltraSparsify} will construct
  a low-stretch spanning tree $T$ of a graph, choose a root vertex $r$,
  apply \texttt{RootedUltraSparsify} to sparsify all edges
  whose path in $T$ contains $r$, and then work
  recursively on the trees obtained by removing the root vertex
  from $T$. 
The root vertex will be chosen to be a tree \textit{splitter},
  where we recall that a vertex $r$ is a \textit{splitter} of a tree
  $T$ if the trees $T^{1}, \dotsc , T^{q}$ obtained by removing $r$
  each have at most half as many vertices as $T$.
The existence of such a vertex was established by Jordan~\cite{Jordan},
  and it is well-known that a tree splitter can be found in linear time.
By making the root a splitter of the tree, we bound the depth of the recursion.
This is both critical for bounding the running time of the algorithm and
  for proving a bound on the quality of the approximation it returns.
For each edge $e$ such that $r \not \in T (e)$, $T (e)$
  is entirely contained in one of $T^{1}, \dotsc , T^{q}$.  
Such edges are sparsified recursively.

\begin{algbox}
\noindent $U = \texttt{UltraSparsify}(G = (V,E), k)$\\
{\bf Condition:} $G$ is connected.

\begin{enumerate}
\item $T = \mathtt{LowStretch} (E)$.
\item Set $t = 517 \cdot 
  \max (1, \log_{2} \eta (E))
  \cdot \ceiling{\log_{2} n} \eta (E) /k$ and 
$p = \left(2 \ceiling{\log \eta (E)} n^{2} \right)^{-1}$.
\item If $t \geq \eta (E)$ then set $A = E-T$;
 otherwise, set
 $A = \mathtt{TreeUltraSparsify} (E - T, t, T, p)$.
\item $U = T \union A$.
\end{enumerate}
\end{algbox}

\begin{algbox}
\noindent $A = \texttt{TreeUltraSparsify}(E', t', T', p)$\\
\begin{enumerate}
\item If $E' = \emptyset$, return $A = \emptyset$.

\item Compute a splitter $r$ of $T'$.

\item Set $E^{r} = \setof{\text{edges $e \in E'$ such that $r \in T'
(e)$}}$ and $t_{r} =\ceiling{t' \eta (E^{r})/\eta (E') }$.

\item If $t_{r} > 1$, set $A^{r} = \mathtt{RootedUltraSparsify}(E^{r}, T', r,
t_{r}, p)$; otherwise, set $A^{r} = \emptyset$.

\item Set $T^{1}, \dotsc , T^{q}$ to be the trees obtained by removing
  $r$ from $T'$.
  Set $V^{1}, \dotsc , V^{q}$ to be the vertex sets of these trees, and
  set $E^{1}, \dotsc , E^{q}$ so that 
  $E^{i} = \setof{(u,v) \in E' : \setof{u,v} \subseteq V^{i}}$.

\item For $i = 1, \dotsc , q$, set
\begin{itemize}
\item [] 
$A = A^{r} \union \texttt{TreeUltraSparsify}(E^{i}, t' \eta (E^{i})/\eta (E') , T^{i},p)$.
\end{itemize}
\end{enumerate}
\end{algbox}

\begin{theorem}[Ultra-Sparsification]\label{thm:ultra}
On input a weighted, connected $n$-vertex graph $G = (V,E)$ and $k \geq 1$,
    \texttt{UltraSparsify}$(E,k)$ returns a set of edges 
    $U = T \union A \subseteq E$ such that
    $T$ is a spanning tree of $G$,
    $U \subseteq E$, and
    with probability at least $1-1/2n$, 
\begin{equation}\label{eqn:ultra1}
U \pleq E \pleq k U,
\end{equation}
and
\begin{equation}\label{eqn:ultra2}
\sizeof{A} \leq O\left(\frac{m}{k} \log^{c_{2}+5} n\right),
\end{equation}
where $m = \sizeof{E}$.
Furthermore, \texttt{UltraSparsify} runs in expected time $O (m\log^{c} n)$, for some constant $c$.
\end{theorem}

We remark that this theorem is very loose when $m/k \geq n$.
In this case, the calls made to \texttt{decompose} by \texttt{RootedUltraSparsify} 
  could have $t \geq n$, in which case \texttt{decompose} will just return singleton
  sets, and the output of \texttt{RootedUltraSparsify} will essentially just
  be the output of \texttt{Sparsify2} on $E^{r}$.
In this case, the upper bound in \eqref{eqn:ultra2} can be very loose.

\begin{proof}
We first dispense with the case $t \geq \eta (E)$.
In this case, \texttt{UltraSparsify} simply returns the graph $E$,
  so \eqref{eqn:ultra1} is trivially satisfied.
The inequality $t \geq \eta (E)$ implies $k \leq O (\log^{2} n)$,
  so \eqref{eqn:ultra2} is trivially satisfied as well.

At the end of the proof, we will use the inequality $t < \eta (E)$.
It will be useful to observe that every time \texttt{TreeUltraSparsify}
  is invoked, 
\[
t' = t \eta (E') / \eta (E).
\]
To apply the analysis of \texttt{RootedUltraSparsify}, we must have
\[
  t^{r} \leq \ceiling{\eta (E^{r})}.
\]
This follows from
\[
  t^{r} = \ceiling{t' \eta (E^{r}) / \eta (E')}
  = \ceiling{t \eta (E^{r}) / \eta (E)}
 \leq  \ceiling{\eta (E^{r})},
\]
as \texttt{TreeUltraSparsify} is only called if $t < \eta (E)$.

Each vertex of $V$ can be a root in a call to \texttt{RootedUltraSparsify}
  at most once, so this subroutine is called at most $n$ times
  during the execution of \texttt{UltraSparsify}.
Thus by Lemma~\ref{lem:rootedUltra}, with probability at least
\[
  1 - n \ceiling{\log_{2} \eta (E)} p = 1 - 1/2n,
\]
  every graph $E_{s}$ returned by a call to \texttt{RootedUltraSparsify}
  satisfies \eqref{eqn:rootedUltraSize} and \eqref{eqn:rootedUltraLeq}.
Accordingly, we will assume both of these conditions hold for the rest of our
  analysis.

We now prove the upper bound on the number of edges in $A$.
During the execution of \texttt{UltraSparsify}, many vertices become the root
  of some tree.
For those vertices $v$ that do not, set $t_{v} = 0$.
By \eqref{eqn:rootedUltraSize},
\begin{equation}\label{eqn:ultraSize1}
\sizeof{A} 
= \sum_{r \in V : t_{r} >1} \sizeof{A^{r}}
 \leq c_{1} \log^{c_{2}} (n/p) \max (1, \ceiling{\log_{2} \eta (E)})
  \sum_{r\in V: t_{r}>1}  t_{r} .
\end{equation}

As $\ceiling{z}\leq 2z$ for $z \geq  1$ and 
  $E^{r_{1}} \cap E^{r_{2}} = \emptyset$
  for each $r_{1}\neq  r_{2}$,
\[
\sum_{r\in V: t_{r}>1} t_{r} 
=
 \sum_{r\in V:t_{r}>1}\ceiling{\frac{\eta (E^{r})}{\eta (E)}t}
\leq 
\sum_{r\in V:t_{r}>1}\frac{2\eta (E^{r})}{\eta (E)}t
\leq  2t.
\]
Thus,
\begin{align*}
\eqref{eqn:ultraSize1} 
& \leq  2 c_{1} \log^{c_{2}} (n/p) \ceiling{\log_{2} \eta (E)} t\\
& \leq 2 c_{1} \log^{c_{2}} (n/p) \ceiling{\log_{2} \eta (E)}
517 \cdot \log_{2} \eta (E)\cdot \ceiling{\log_{2} n} \eta (E) /k\\
& \leq  O \left( \frac{m}{k} \log^{c_{2} + 5} n \right),
\end{align*}
where the last inequality uses $\eta (E) = O (m \log n \ \log \log n) = 
  O (m \log^{2} n)$
  from Theorem~\ref{thm:lowStretch} and $\log m = O (\log n)$.

We now establish \eqref{eqn:ultra1}.
For every vertex $r$ that is ever 
  selected as a tree splitter in line 2 of \texttt{TreeUltraSparsify},
  let $T^{r}$ be the tree $T'$ of which $r$ is a splitter,
  and let $E^{r}$ denote the set of edges and
  $t_{r}$ be the parameter set in line 3.
Observe that $\union_{r} E^{r} = E - T$.
Let
\[
  \beta_{r} = 4 \eta (E^{r}) / t_{r},
\]
and note this is the parameter used in the analysis
  of \texttt{RootedUltraSparsify} in Lemma~\ref{lem:rootedUltra}.
If $t_{r} > 1$, let $A^{r}$ be the set of edges returned
  by the call to \texttt{RootedUltraSparsify}.
By Lemma~\ref{lem:rootedUltra}, \texttt{RootedUltraSparsify} 
  returns a set of edges $A^{r}$ satisfying
\begin{equation}\label{eqn:TUSresult}
E^{r} \pleq 
  \left(3 \beta_{r}  + 126 \beta_{r}  \max (1, \log_{2} \eta (E^{r}))\right)
  \cdot T^{r} 
  + 120 \beta_{r}  \cdot A^{r}.
\end{equation}
On the other hand, if $t_{r} = 1$ and so $A^{r} = \emptyset$,
  then $\beta_{r} = 4 \eta (E^{r})$.
We know that \eqref{eqn:TUSresult} is satisfied in this case because
  $E^{r} \pleq \eta (E^{r}) T^{r}$ (by \eqref{eqn:treeStretch}).
If $t_{r} = 0$, then $E^{r} = \emptyset$ and \eqref{eqn:TUSresult} is
  trivially satisfied.
As $t_{r} = \ceiling{t \eta (E^{r}) / \eta (E)},$
\[
  \beta_{r} \leq 4 \eta (E) / t.
\]
We conclude
\[
E^{r} \pleq 
  129 \beta_{r} \max (1, \log_{2} \eta (E^{r}))
  \cdot T^{r} 
  + 120 \beta_{r}  \cdot A^{r}
\pleq 
  516 (\eta (E)/ t) \max (1,\log_{2} \eta (E^{r})) T^{r} + 120 (\eta (E)/t) A^{r}.
\]
Adding $T$, 
  summing over all $r$, and remembering $\eta (E^{r}) \leq \eta (E)$, we obtain
\[
T + (E-T) \pleq T + 516 (\eta (E)/t) \max (1,\log_{2} \eta (E))
\sum_{r} T^{r}
    + 120 (\eta (E)/t) A.
\]
As $r$ is always chosen to be a splitter of the tree input to
  \texttt{TreeUltraSparsify},
  the depth of the recursion is at most $\ceiling{\log_{2} n}$.
Thus, no edge of $T$ appears more than $\ceiling{\log_{2} n}$
  times in the sum $\sum_{r} T^{r}$, and we may conclude
\begin{align*}
T + (E-T) &
  \pleq 
 T + 516 (\eta (E)/t) \max (1,\log_{2} \eta (E)) \ceiling{\log_{2} n} T  + 120 (\eta (E)/t) A
\\
&  \pleq 
 517 (\eta (E)/t) \max (1,\log_{2} \eta (E)) \ceiling{\log_{2} n} T  + 120 (\eta (E)/t) A
\\
& \pleq 
k (T + A)
\\
& =
k U,
\end{align*}
where the second inequality follows from $t \leq \eta (E)$,
and the third inequality follows from the value chosen for $t$
  in line 2 of \texttt{UltraSparsify}.

To bound the expected running time of \texttt{UltraSparsify},
 first observe that the call to \texttt{LowStretch} takes time $O (m \log n \ \log \log n)$.
Then, note that
 the routine \texttt{TreeUltraSparsify} is recursive, the recursion
  has depth at most $O (\log n)$, and
  all the graphs being processed by \texttt{TreeUltraSparsify}
  at any level of the recursion are disjoint.
The running time of \texttt{TreeUltraSparsify} is dominated by the
  calls made to \texttt{Sparsify2} inside
  \texttt{RootedUltraSparsify}.
Each of these takes nearly-linear expected time, so the overall
  expected running time of \texttt{TreeUltraSparsify} is $O (m \log^{c} n) $, for some constant $c$.
\end{proof}

\bibliographystyle{alpha}
\bibliography{precon}

\appendix

\section{Gremban's reduction}\label{sec:gremban}

Gremban~\cite{Gremban} (see also~\cite{MaggsEtAl}) provides
  the following method for handling
  positive off-diagonal entries.
If $A$ is a \sdd-matrix,
  then Gremban decomposes
  $A$ into $D + A_{n} + A_{p}$, where
  $D$ is the diagonal of $A$,
  $A_{n}$ is the matrix containing all the
  negative off-diagonal entries of $A$,
  and $A_{p}$ contains all the positive off-diagonals.
Gremban then considers the linear system
\[
\Ahat \left(
\begin{array}{l}
\xx_{1}\\
\xx_{2}
\end{array}
 \right)
= 
\bbhat,
\qquad \text{where} \qquad 
\Ahat =   \left[
\begin{array}{ll}
  D + A_{n} & -A_{p}\\
  -A_{p} & D + A_{n}
\end{array}
 \right]
\quad \text{and} \quad
\bbhat = \left(
\begin{array}{l}
\bb\\
-\bb
\end{array}
 \right),
\]
and observes that 
  $\xx = (\xx_{1} - \xx_{2})/2$
  will be the solution to
  $A \xx = \bb $, if a solution exists.
Moreover, approximate solutions of Gremban's system yield approximate
  solutions of the original:

\[
\norm{
\left(
\begin{array}{l}
\xx_{1}\\
\xx_{2}
\end{array}
 \right)
-
\pinv{\Ahat}
\bbhat
}
\leq 
\epsilon
\norm{
\pinv{\Ahat}
\bbhat
}
\qquad \text{implies} \qquad 
  \norm{\xx - \pinv{A} \bb }
\leq
  \epsilon \norm{\pinv{A} \bb},
\]
where again $\xx = (\xx_{1} - \xx_{2})/2$.
Thus we may reduce the problem of solving a linear system in
  a \sdd-matrix into that of solving a linear system in a
  \sddm-matrix that is at most twice as large and has at most twice
  as many nonzero entries.

\section{Computing the stretch}\label{sec:stretch}

We now show that given a weighted graph $G = (V,E)$ 
  and a spanning tree $T$ of $G$,
  we can compute $\stretch{T}{e}$ for every edge 
  $e\in E$ in $O ((m + n)\log  n)$ time, where
  $m=\sizeof{E}$ and $n=\sizeof{V}$.

For each pair of vertices $u,v \in V$, let $\res{u,v}$ be 
   the resistance of $T (u,v)$, the path in $T$ connecting $u$ and $v$.
We first observe that for an arbitrary $r\in V$, we can compute
  $\res{v,r}$ for all $v\in V$ in $O (n)$ time by
  a top-down traversal on the rooted tree obtained from 
  $T$ with root $r$.
Using this information, we can compute the stretch of all
  edges in $E_{r} = \setof{\text{edges $e \in E$ such that $r \in T
(e)$}}$ in time $O (\sizeof{E_{r}})$.

We can then use tree splitters in the same manner as in \texttt{TreeUltraSparsify}
  to compute the stretch of all edges in $E$ in $O ((m+n)\log n)$ time.
That is, in linear time we can identify a vertex $r$ such that the removal of $r$
  from of the tree produces sub-trees each of which have at most half as many
  vertices as the original.
We then treat $r$ as the root and compute the stretch of all edges
  in $E_{r}$.
We note that these are exactly the edges
 whose endpoints
  are in different components of the forest obtained by removing all edges attached
  to vertex $r$.
We then recurse on the sub-trees obtained by removing the vertex $r$ from the tree.
This algorithm performs a linear amount of work on each level of the recursion,
  and our choice of splitters as roots of the trees guarantees that there
  are at most $\ceiling{\log_{2} n}$ levels of recursion.

\section{Decomposing Trees}\label{sec:decompose}

The pseudo-code for \texttt{decompose} appears on the next page.
The algorithm performs a depth-first traversal of the tree, greedily
  forming sets $W_{j}$ once they are attached to edges whose sum of $\eta$ values
  exceeds a threshold $\phi$.
The traversal is performed by the subroutine \texttt{sub}.
This routine is first called from the root of the tree, and it
  then recursively calls itself on the children of its input vertex
  before processing its input vertex.
The routine \texttt{sub} returns a set of vertices, $U$, along with
  $F$, 
  the set of edges touching vertices in $U$, and $w$, the sum
  of $\eta$ over the edges in $F$.
The vertices in $U$ are those in the sub-tree rooted at $v$
  which \texttt{sub} did not place into a set $W_{j}$.
While \texttt{sub} is gathering sets of vertices $U_{sub}$, the edges they are attached to
  are stored in $F_{sub}$, and the sum of the value of $\eta$ on these edges
  is stored in $w_{sub}$.

There are four lines on which the routine \texttt{sub} can form a set $W_{j}$:
  on lines 3.c.ii, 6.b, 7.b or 7.e. 
When the routine \texttt{sub} forms a set $W_{j}$ of vertices, those vertices
  form a sub-tree of $T$.
If the set is formed on line 7.e, then this subtree will be a singleton consisting of only one vertex.
All of the edges attached to those vertices, except possibly for the root
  of that sub-tree, are assigned to $W_{j}$ by $\rho$.
If the set $W_{j}$ is formed in line 6.b or 7.e, then all of the edges attached to
  the root are assigned to $W_{j}$.
Otherwise, the edges attached to the root are not assigned to $W_{j}$ by $\rho$,
  unless they happen to also be attached to another vertex in $W_{j}$.

\begin{figure}[htop]
\begin{algbox}
\noindent $(\setof{W_{1}, \dots , W_{h}}, \rho ) = 
    \texttt{decompose}(T, E, \eta, t)$\\

{\bf Comment}: $h$, $\rho $, and the $W_{i}$'s are treated as global variables.
\begin{enumerate}
\setlength{\itemsep}{0pt}
\setlength{\topsep}{0pt}
\setlength{\parsep}{0pt}
\setlength{\parskip}{0pt}
\setlength{\partopsep}{0pt}
\item Set $h = 0$.  \qquad  \qquad \qquad \qquad \qquad  \textit{($h$ is incremented as sets are created)}
\item For all $e\in E$, set $\rho (e) = \emptyset$.
\item Set $\phi = 2 \sum_{e \in E} \eta (e) / t$. \qquad \textit{(the threshold above which sets are formed)}
\item $(F, w, U) = \texttt{sub}(r)$.
\item If $U\neq \emptyset$,
   \begin{enumerate}
   \setlength{\itemsep}{0pt}
   \setlength{\topsep}{0pt}
   \setlength{\parsep}{0pt}
   \setlength{\parskip}{0pt}
   \setlength{\partopsep}{0pt}
   \item $h = h + 1.$
   \item $W_{h} = U$.
   \item For all $e \in F$, set $\rho (e) = \rho (e) \union \setof{W_{h}}$.
   \end{enumerate}
\end{enumerate}

\vskip 20pt

\noindent $(F, w, U) = \texttt{sub}(v)$

\qquad \textit{$U$ is a set of vertices, $F$ is the set of edges attached to $U$, and $w$ is the sum of $\eta$ over $F$}
\begin{enumerate}
\setlength{\itemsep}{0pt}
\setlength{\topsep}{0pt}
\setlength{\parsep}{0pt}
\setlength{\parskip}{0pt}
\setlength{\partopsep}{0pt}
\item Let $v_{1}, \dots , v_{s}$ be the children of $v$.
\item Set $w_{sub} = 0$, $F_{sub} = \emptyset$ and $U_{sub} = \emptyset$.
\item For $i = 1, \dots , s$
  \begin{enumerate}
\setlength{\itemsep}{0pt}
\setlength{\topsep}{0pt}
\setlength{\parsep}{0pt}
\setlength{\parskip}{0pt}
\setlength{\partopsep}{0pt}
  \item $(F_{i}, w_{i}, U_{i}) = \texttt{sub} (v_{i})$.
  \item $w_{sub} = w_{sub} + w_{i}$, $F_{sub} = F_{sub} \union F_{i}$, $U_{sub} = U_{sub} \union U_{i}$.
  \item If $w_{sub} \geq  \phi $,
    \begin{enumerate}
\setlength{\itemsep}{0pt}
\setlength{\topsep}{0pt}
\setlength{\parsep}{0pt}
\setlength{\parskip}{0pt}
\setlength{\partopsep}{0pt}
    \item $h = h + 1$.
    \item Set $W_{h} = U_{sub} \union \setof{v}$.
    \item For all $e \in F_{sub}$, set $\rho (e) = \rho (e) \union \setof{W_{h}}$.
    \item Set $w_{sub} = 0$, $F_{sub} = \emptyset$ and $U_{sub} = \emptyset$.
    \end{enumerate}
  \end{enumerate}
  
  \item Set $F_{v} = \setof{(u,v) \in E}$, the edges attached to $v$.
  \item Set $w_{v} = \sum_{e \in F_{v}} \eta (e)$.
  \item If $\phi \leq w_{v} + w_{sub} \leq  2 \phi $,
    \begin{enumerate}
\setlength{\itemsep}{0pt}
\setlength{\topsep}{0pt}
\setlength{\parsep}{0pt}
\setlength{\parskip}{0pt}
\setlength{\partopsep}{0pt}
    \item $h = h + 1$.
    \item Set $W_{h} = U_{sub} \union \setof{v}$.
    \item For all $e \in F_{sub} \union F_{v}$, set $\rho (e) = \rho (e) \union \setof{W_{h}}$.
    \item Return $(\emptyset, 0, \emptyset)$.
    \end{enumerate}
  \item If $w_{v} + w_{sub} > 2 \phi$,
    \begin{enumerate}
\setlength{\itemsep}{0pt}
\setlength{\topsep}{0pt}
\setlength{\parsep}{0pt}
\setlength{\parskip}{0pt}
\setlength{\partopsep}{0pt}
    \item $h = h + 1$.
    \item Set $W_{h} = U_{sub} \union \setof{v}$.  
    \item For all $e \in F_{sub}$, set $\rho (e) = \rho (e) \union \setof{W_{h}}$.
    \item $h = h + 1$.
    \item Set $W_{h} = \setof{v} $. \qquad \qquad \qquad \qquad \qquad \textit{(create a singleton set)}
    \item For all $e \in F_{v}$, set $\rho (e) = \rho (e) \union \setof{W_{h}}$.
    \item Return $(\emptyset, 0, \emptyset)$.
  \end{enumerate}
  \item Return $(F_{sub} \union F_{v}, w_{sub} + w_{v}, U_{sub} \union \setof{v})$
\end{enumerate}
\end{algbox}
\end{figure}

There are three ways that \texttt{sub} could decide to create sets $W_{j}$.
The first is if some subset of the children of $v$ return sets $F_{i}$
  whose values under $\eta$ sum to more than $\phi$.
In this case, \texttt{sub} collects the corresponding sets of vertices,
  along with $v$, into a set $W_{j}$.
But, the edges attached to $v$ do not necessarily get assigned by $\rho$ to $W_{j}$
  as $v$ was merely included in the set to make it a connected subtree.
The second is if the sum of $\eta$ over the edges attached to $v$ and the sets
  $U_{i}$ returned by a bunch of the children exceeds $\phi$, but is less than $2 \phi$.
In this case, $v$ and all those sets $U_{i}$ are bundled together into a set $W_{j}$
  in line 6.b.
Finally, if the sum of $\eta$ over the edges attached to $v$ and the 
 the sets
  $U_{i}$ returned by a bunch of the children exceeds $ 2\phi$, then $v$
  is added as a singleton set, and a set is created containing the union of
  $v$ with those sets $U_{i}$.

We assume that some vertex $r$ has been chosen to be the root of the tree.
This choice is used to determine which nodes in the tree are children
  of each other.

\begin{proof}[{\rm\bf Proof} of Theorem~\ref{thm:decomposeTree}]
As algorithm \texttt{decompose} traverses the tree $T$
   once and visits each edge in $E$ once, it runs in linear time.

In our proof, we will say that an edge $e$ is assigned to a set $W_{j}$
  if $W_{j} \in \rho (e)$.
To prove part $(a)$ of the theorem, we use the following observations:
  If $W_{j}$ is formed in step 3.c.ii or step 6.b, then
  the sum of $\eta$ over edges assigned to $W_{j}$ is at least $\phi $,
  and if $W_{j}$ is formed in step 7.b, then the sum of $\eta$
  of edges incident to $W_{j}$ and $W_{j+1}$ (which is a singleton)
  is at least $2\phi $.
Finally, if a set $W_{h}$ is formed in line 5.b of \texttt{decompose},
  then the sum of $\eta $ over edges assigned to $W_{h}$ is greater than zero.
But, at most one set is formed this way.
As each edge is assigned to at most two sets in $W_{1},\dots ,W_{h}$,
we may conclude
\[
  2 \sum_{e \in E} \eta (e) > (h-1) \phi ,
\]
which implies $t > h-1$.
As both $t$ and $h$ are integers, this implies $t \geq h$.

We now prove part $(b)$.
First, observe that steps 6 and 7 guarantee that
  when a call to $\mathtt{sub} (v)$ returns a triple
  $(F, w, U)$, 
\[
w = \sum_{e \in F} \eta (e) < \phi .
\]
Thus, when a set $W_{h}$ is formed in step 3.c.ii, we know that
  the sum of $\eta$ over edges assigned to $W_{h}$
  equals $w_{sub}$ and is at most $2 \phi$.
Similarly, we may reason that $w_{sub} < \phi$ at step 4.
If a set $W_{h}$ is formed in step 6.b, the sum of $\eta$
  over edges associated with $W_{h}$ is $w_{v} + w_{sub}$,
  and must be at most $2 \phi$.
If a set $W_{h}$ is formed in step 7.b, the sum of $\eta$
  over edges associated with $W_{h}$ is $w_{sub}$,
  which we established is at most $\phi$.
As the set formed in step 7.e is a singleton, we do not need to
  bound the sum of $\eta$ over its associated edges.
\end{proof}

\begin{lemma}\label{lem:planar}
Suppose $G= (V,E)$ is a planar graph, $\pi$ is a planar embedding
  of $G$, $T$ is a spanning tree of $G$, and $t > 1$ is an integer.
Let $(\setof{W_{1},\dotsc ,W_{h}},\rho ) = \texttt{decompose}
  (T,E,\eta,t)$ with the assumption that in Step 1 of \texttt{sub}, 
  the children $v_{1},\dotsc ,v_{s}$ of $v$ always appear in clock-wise order
  according to $\pi$.
Then the graph
   $G_{\setof{W_{1},\dotsc ,W_{h}}} = \left(\setof{1,\dotsc ,h}, \setof{(i,j): \exists\  e\in E, \rho (e) = \setof{W_{i},W_{j}}}\right)$
is planar.
\end{lemma}
\begin{proof}
Recall that the contraction of an edge $e = (u,v)$ in a planar
  graph $G = (V,E)$
  defines a new graph $\left(V-\setof{u}, E \cup \setof{(x,v) :
(x,u)\in E} -   \setof{(x,u) \in E} \right)$.
Also recall that 
  edge deletions and edge contractions preserve  planarity.

We first prove the lemma in the special case in which
  the sets  $W_{1},\dotsc ,W_{h}$ are disjoint.
For each $j$, let $T_{j}$ be the graph induced on $T$ by $W_{j}$.
As each $T_{j}$ is connected,
  $G_{\setof{W_{1},\dotsc ,W_{h}}}$ is a subgraph of the graph
  obtained by contracting all the edges in each subgraph $T_{j}$.
Thus in this special case $G_{\setof{W_{1},\dotsc ,W_{h}}}$ is planar.

We now analyze the general case, recalling that the sets $W_{1}, \dotsc , W_{h}$
  can overlap.
However, the only way sets $W_{j}$ and $W_{k}$ with $j < k$ 
  can overlap is if the set $W_{j}$ was formed at step 3.c.ii, and the
  vertex $v$ becomes part of $W_{k}$ after it is returned by a call to \texttt{sub}.
In this situation, no edge is assigned to $W_{j}$ for having $v$ as an end-point.
That is, the only edges of form $(x,v)$ that can be assigned to $W_{j}$
  must have $x \in W_{j}$.
So, these edges will not appear in $G_{\setof{W_{1},\dotsc ,W_{h}}}$.

Accordingly, for each $j$ we define
\[
X_{j} = \begin{cases}
W_{j} - v  & \text{if $W_{j}$ was formed at step 3.c.ii, and}
\\
W_{j}  & \text{otherwise}.
\end{cases}
\]
We have shown that $G_{\setof{W_{1},\dotsc ,W_{h}}} = G_{\setof{X_{1},\dotsc ,X_{h}}}$.
Moreover, the sets $X_{1}, \dotsc , X_{h}$ are disjoint.
Our proof would now be finished, if only each subgraph of $G$ induced
  by a set $X_{j}$ were connected.
While this is not necessarily the case, we can make it the case by adding
  edges to $E$.

The only way the subgraph of $G$ induced on a set $X_{j}$ can fail to
  be connected is if $W_{j}$ is formed at line 3.c.ii from 
   the union of $v$ with a collection
  sets $U_{i}$ for $i_{0} \leq i \leq i_{1}$ returned by recursive calls
  to \texttt{sub}.
Now, consider what happens if we add edges of the form $(v_{i}, v_{i+1})$	 
  to the graph for $i_{0} \leq i < i_{1}$, whenever they are not already present.
As the vertices $v_{i_{0}}, \dotsc , v_{i_{1}}$ appear in clock-wise order around $v$,
  the addition of these edges preserves the planarity of the graph.
Moreover, their addition makes the induced subgraphs on each set $X_{j}$ connected,
  so we may conclude that $G_{\setof{X_{1},\dotsc ,X_{h}}}$ is in fact planar.
\end{proof}

\end{document}